\documentclass[12pt]{amsart}
\usepackage{latexsym, amsmath, amscd, amssymb, amsthm}
\usepackage{bm, mathrsfs, euscript}
\usepackage{subfigure, epsfig, psfrag}
\usepackage{pstricks,pst-node,pst-text,pst-tree}
\usepackage[l>=1.2em]{diagrams}
\usepackage[colorlinks=true,linkcolor=blue,citecolor=red]{hyperref}
\usepackage[alphabetic]{amsrefs}
\usepackage[centering, includeheadfoot, hmargin=1.2in, vmargin=0.8in,
  headheight=30.4pt]{geometry}
\newtheorem{lemma}{Lemma}[section]
\newtheorem{theorem}[lemma]{Theorem}
\newtheorem{corollary}[lemma]{Corollary}
\newtheorem{proposition}[lemma]{Proposition}
\theoremstyle{definition}

\newtheorem{example}[lemma]{Example}
\newtheorem{remark}[lemma]{Remark}
\newtheorem{conventions}[lemma]{Conventions}
\newtheorem*{acknowledgements}{Acknowledgements}
\theoremstyle{remark}
\newtheorem*{proof1}{$(1)$~Choosing the line bundles $L_{r-1}$ and
  $L_r$}
\newtheorem*{proof2}{$(2)$~Proof that $I_Q = \bigl( I_R : (\prod_{a
    \in Q_1} y_a)^\infty \bigr)$}
\newtheorem*{proof3}{$(3)$~Proof that $(I_R : B_Y^\infty) = \bigl( I_R
  : (\prod_{a \in Q_1} y_a)^\infty \bigr)$}

\numberwithin{equation}{section}
\renewcommand{\theequation}%
{\arabic{section}.\arabic{lemma}.\arabic{equation}}
  
\newarrow{Equal}{=}{=}{=}{=}{=}
\newarrow{Dashto}{dash}{dash}{dash}{dash}{>}
\renewcommand{\div}{\operatorname{div}}
\newcommand{\kk}{\ensuremath{\Bbbk}} 
\newcommand{\git}{\ensuremath{/\!\!/_{\!\theta}\,}}

\renewcommand{\AA}{\ensuremath{\mathbb{A}}} 
 
\newcommand{\FF}{\ensuremath{\mathbb{F}}} 
\newcommand{\NN}{\ensuremath{\mathbb{N}}} 
\newcommand{\PP}{\ensuremath{\mathbb{P}}} 
\newcommand{\QQ}{\ensuremath{\mathbb{Q}}} 
\newcommand{\RR}{\ensuremath{\mathbb{R}}} 
\newcommand{\VV}{\ensuremath{\mathbb{V}}} 
\newcommand{\ZZ}{\ensuremath{\mathbb{Z}}} 
\DeclareMathOperator{\Amp}{Amp}
\DeclareMathOperator{\CDiv}{CDiv}
\DeclareMathOperator{\Cl}{Cl}
\DeclareMathOperator{\Cir}{Cir}
\DeclareMathOperator{\End}{End}
\DeclareMathOperator{\eval}{ev}

\DeclareMathOperator{\head}{hd}
\DeclareMathOperator{\Hom}{Hom}
\DeclareMathOperator{\inc}{inc}
\DeclareMathOperator{\Wt}{Wt}
\DeclareMathOperator{\mult}{mult}

\DeclareMathOperator{\Ker}{Ker}
\DeclareMathOperator{\Nef}{Nef}

\DeclareMathOperator{\Pic}{Pic}
\DeclareMathOperator{\pic}{pic}
\DeclareMathOperator{\Proj}{Proj}
\DeclareMathOperator{\rad}{rad}

\DeclareMathOperator{\Spec}{Spec}
\DeclareMathOperator{\supp}{supp}
\DeclareMathOperator{\tail}{tl}
\DeclareMathOperator{\modules}{-mod}

\begin{document}

\title[Toric varieties as fine moduli]{Projective Toric Varieties as
  \\ fine moduli spaces of Quiver Representations}

\author[A.~Craw]{Alastair Craw} 
\address{Department of Mathematics\\ University of Glasgow\\
  University Gardens\\ Glasgow\\ G12 8QW\\ United Kingdom}
\email{craw@maths.gla.ac.uk}

\author[G.G.~Smith]{Gregory G. Smith} 
\address{Department of Mathematics and Statistics \\ Queen's
  University \\ Kingston \\ ON \\ K7L 3N6\\ Canada}
\email{ggsmith@mast.queensu.ca}

%\date{17 March 2007}

\begin{abstract}
  This paper proves that every projective toric variety is the fine
  moduli space for stable representations of an appropriate bound
  quiver.  To accomplish this, we study the quiver $Q$ with relations
  $R$ corresponding to the finite-dimensional algebra $\End\bigl(
  \textstyle\bigoplus\nolimits_{i=0}^{r} L_i \bigr)$ where
  $\mathcal{L} := (\mathscr{O}_X,L_1, \dotsc, L_r)$ is a list of line
  bundles on a projective toric variety $X$.  The quiver $Q$ defines a
  smooth projective toric variety, called the multilinear series
  $|\mathcal{L}|$, and a map $X \rTo |\mathcal{L}|$.  We provide
  necessary and sufficient conditions for the induced map to be a
  closed embedding.  As a consequence, we obtain a new geometric
  quotient construction of projective toric varieties.  Under slightly
  stronger hypotheses on $\mathcal{L}$, the closed embedding
  identifies $X$ with the fine moduli space of stable representations
  for the bound quiver $(Q,R)$.
\end{abstract}

\maketitle

\vspace*{-1em}

\section{Introduction} 
\label{sec:introduction}

The dictionary between the geometry of a moduli space $X$ and the
family of objects classified by $X$ lies at the heart of modern
algebraic geometry.  Fine moduli spaces, although substantially rarer
than coarse ones, provide the fundamental example of this
correspondence.  A scheme $X$ is a fine moduli space for the
equivalence classes of some objects if and only if there is a
universal family of the selected objects over $X$ such that every
other family of objects is induced from the universal one by a unique
morphism to $X$.  The hyperplane bundle on $\PP^d$ and the
tautological vector bundle on a Grassmannian are the classic examples
of universal families.  A universal family is a powerful tool for
studying the geometry of $X$ as illustrated by the forgetful morphism
between moduli spaces of pointed stable curves, or the Fourier-Mukai
transform on abelian varieties.  The primary goal of this paper is to
realize every projective toric variety as a fine moduli space of
stable representations for an appropriate bound quiver.

To be more precise, let $X$ be a projective toric variety over a field
$\kk$ and consider a list $\mathcal{L} := (\mathscr{O}_X, L_1, \dotsc,
L_r)$ of line bundles on $X$; for convenience set $L_0 :=
\mathscr{O}_X$.  The finite-dimensional $\kk$-algebra $\End \bigl(
\textstyle\bigoplus\nolimits_{i=0}^{r} L_i \bigr)$ is encoded by a
quiver $Q$, called the complete quiver of sections for $\mathcal{L}$,
together with an ideal of relations $R$ in the path algebra $\kk Q$.
We associate to the quiver $Q$ a unimodular, projective toric variety
$|\mathcal{L}|$ called the multilinear series of $Q$.  The variety
$|\mathcal{L}|$ can be defined combinatorially, by geometric invariant
theory, or via representation theory; see
Proposition~\ref{pro:smooth}.  Since the isomorphism $\kk Q/R \cong
\End \bigl( \bigoplus_{i=0}^{r} L_i \bigr)$ identifies arrows in $Q$
with global sections of line bundles, the quiver $Q$ induces a map
from $X$ to $|\mathcal{L}|$.  We prove the following:

\begin{theorem}
  \label{thm:1}
  If $L_1, \dotsc, L_r$ are basepoint-free line bundles on $X$, then
  the induced map $\varphi_{|\mathcal{L}|} \colon X \rTo |\mathcal{L}|$
  is a morphism and the image is presented as a geometric quotient.
\end{theorem}

\noindent
When $r = 1$, $\varphi_{|\mathcal{L}|}$ coincides with the morphism
from $X$ to the linear series $|L_1|$.

Lists $\mathcal{L}$ for which the induced morphism $\varphi_{|
  \mathcal{L}|}$ is a closed embedding are ubiquitous; see
Proposition~\ref{pro:veryample}.  Hence, Theorem~\ref{thm:1} produces
a wealth of new geometric quotient constructions for a projective
toric variety $X$.  In particular, these \emph{geometric} quotients
provide new ``homogeneous coordinates'' for the points on $X$.  In
contrast with \cites{Cox, Kajiwara}, these quotient constructions are
not intrinsic to the toric variety; they depend on the choice of line
bundles in $\mathcal{L}$.  Although we recover the quotient
constructions in \cites{Cox, Kajiwara} for some toric varieties and
particular lists $\mathcal{L}$, the homogeneous coordinate systems
arising from a multilinear series $|\mathcal{L}|$ are typically
larger.  Examples suggest that some of these larger coordinate systems
appear naturally in the quantum cohomology of $X$ and the derived
category of coherent sheaves on $X$.

To achieve the primary goal, we relate the image
$\varphi_{|\mathcal{L}|}(X)$ to an important subscheme of
$|\mathcal{L}|$.  From the viewpoint of representation theory, the
multilinear series $|\mathcal{L}|$ is the fine moduli space of
$\vartheta$-stable representations with dimension vector $(1, \dotsc,
1)$ for the quiver $Q$ where $\vartheta$ is a distinguished weight on
$Q$; see Proposition~\ref{pro:smooth}.  The ideal of relations $R$ in
the path algebra $\kk Q$ determines a subscheme of $|\mathcal{L}|$
that coincides with the fine moduli space $\mathcal{M}_\vartheta(Q,R)$
of $\vartheta$-stable representations with dimension vector $(1,
\dotsc, 1)$ for the bound quiver $(Q,R)$.  In other words, the variety
$\mathcal{M}_\vartheta(Q,R)$ classifies certain finite-dimensional
modules over the $\kk$-algebra $\End \bigl( \bigoplus_{i=0}^{r} L_i
\bigr)$.  The universal $\vartheta$-stable representation of $(Q,R)$
over $\mathcal{M}_\vartheta(Q,R)$ decomposes into a direct sum of line
bundles called the tautological line bundles.  Our main results can be
summarized as follows:

\begin{theorem}
  \label{thm:2}
  Let $X$ be a projective toric variety.  There exist (many) lists
  $\mathcal{L}$ of line bundles on $X$ such that the induced morphism
  $\varphi_{|\mathcal{L}|} \colon X \rTo |\mathcal{L}|$ identifies $X$
  with the fine moduli space $\mathcal{M}_\vartheta(Q,R)$.  Moreover,
  the tautological line bundles on $\mathcal{M}_\vartheta(Q,R)$
  coincide with the line bundles in the list $\mathcal{L}$.
\end{theorem}
  
\noindent
This fine moduli interpretation yields a functorial approach to
projective toric varieties.  Specifically, it allows one to describe
the data needed to specify a map from a scheme to a projective toric
variety as in \cites{Cox2,Kajiwara}.

Theorem~\ref{thm:2} also helps clarify the relationship between
descriptions of the derived category $D^{b}(\mathcal{O}_X\modules)$
and realizations of $X$ as a fine moduli space of quiver
representations.  \cite{Bondal} shows that
$D^{b}(\mathcal{O}_X\modules)$ is equivalent to the derived category
of finite-dimensional modules over $\End \bigl( \bigoplus_{i=0}^{r}
\mathscr{F}_i \bigr)$ if and only if the coherent sheaves
$\mathscr{F}_{i}$ form a complete strong exceptional collection on
$X$.  On certain toric quiver varieties, \cite{AH} describe such
collections in which the $\mathscr{F}_i$ are line bundles; toric
quiver varieties are fine moduli spaces of quiver representations.
The influential \cite{King2} constructs complete strong exceptional
collections of line bundles on several smooth toric surfaces by
realizing the surfaces as fine moduli spaces of stable representations
of a bound quiver.  Given a smooth projective variety with a complete
strong exceptional collection of line bundles, \cite{BergmanProudfoot}
establishes that the variety is isomorphic to a connected component of
a corresponding moduli space of stable quiver representations.
\cite{Kawamata} proves that every toric variety has a complete
exceptional collection of coherent sheaves and, in contrast,
\cite{HillePerling} exhibits a smooth toric surface that does not have
a complete strong exceptional collection of line bundles.  In this
context, Theorem~\ref{thm:2} clearly differentiates between a fine
moduli interpretation of a variety and the existence of a complete
strong exceptional collection of line bundles.
 
This paper is organized as follows.  Our notation and some standard
results from toric geometry and quiver theory are described in
\S\ref{sec:background}.  In \S\ref{sec:quiver}, we define a quiver of
sections and its associated multilinear series.  This generalizes the
classical notion of a linear series from a single line bundle to a
list of line bundles.  The induced map to the multilinear series is
studied in \S\ref{sec:multilinear}.  In particular, we give necessary
and sufficient conditions for the induced map to be a morphism or a
closed embedding.  Finally, \S\ref{sec:boundreps} examines
representations of a bound quiver of sections and establishes our main
results.

\begin{acknowledgements}
  We thank Diane Maclagan, Sorin Popescu, Nicholas Proudfoot and Karen
  Smith for stimulating conversations.  We also thank an anonymous
  referee for valuable comments and suggestions.  This paper owes much
  to experiments made with \cite{M2}.  The second author is partially
  supported by NSERC.
\end{acknowledgements}

\section{Background and Notation} 
\label{sec:background}

We collect here standard definitions, results and notation.  In this
paper, $\NN$ denotes the nonnegative integers and $\kk$ is an
algebraically closed field of characteristic zero.

\subsection{Toric Varieties}
\label{sub:toric}

Let $X$ be a projective toric variety over $\kk$ determined by a
strongly convex rational polyhedral fan $\Sigma_{X} \subseteq N_X
\otimes_{\ZZ} \RR \cong \RR^{d}$ where $N_X$ is a lattice of rank $d$.
The dual lattice is $M_X := \Hom_{\ZZ}(N_X,\ZZ)$ and $T_X := N_X
\otimes_{\ZZ} \kk^{*}$ is the algebraic torus acting on $X$.  The
$i$-dimensional cones of $\Sigma_{X}$ form the set $\Sigma_{X}(i)$.
Since $X$ is projective, $\bigcup_{\sigma \in \Sigma_{X}} \sigma =
N_{X} \otimes_{\ZZ} \RR$ and $\Sigma_{X}(d)$ is the set of maximal
cones.

Each $\rho \in \Sigma_{X}(1)$ corresponds to an irreducible
$T_X$-invariant Weil divisor $D_{\rho}$ on $X$.  These divisors
generate the free abelian group $\ZZ^{\Sigma_{X}(1)}$ of
$T_X$-invariant Weil divisors and the semigroup $\NN^{\Sigma_{X}(1)}$
of effective $T_X$-invariant Weil divisors.  The quotient of
$\ZZ^{\Sigma_{X}(1)}$ by the subgroup of principal divisors is the
class group (or Chow group) $\Cl(X)$.  The $T_X$-invariant Cartier
divisors $\CDiv(X)$ form a subgroup of $\ZZ^{\Sigma_{X}(1)}$.
Moreover, there is a commutative diagram
\begin{equation} 
  \label{eq:fundiagram}
  \begin{diagram}[size=2em]
    0 & \rTo & M_X & \rTo & \CDiv(X) & \rTo & \Pic(X) & \rTo
    & 0 \\
    & & \dEqual & & \dTo & & \dTo & & \\
    0 & \rTo & M_X & \rTo & \ZZ^{\Sigma_{X}(1)} & \rTo & \Cl(X) & \rTo
    & 0
  \end{diagram}
\end{equation}
where the rows are exact and the vertical arrows are inclusions; see
\S3.4 in \cite{Fulton}.  The projection from $\ZZ^{\Sigma_{X}(1)}$ to
$\Cl(X)$ is denoted by $u \mapsto [u]$ and the inclusion of $\Pic(X)$
into $\Cl(X)$ is also denoted by $L \mapsto [L]$.  For a line bundle
$L$ on $X$ and a global section $s \in H^0(X,L)$, $\div(s)$ denotes
the effective Cartier divisor determined by $s$.

The total coordinate ring of $X$ is the polynomial ring $S_X :=
\kk[x_{\rho} : \rho \in \Sigma_{X}(1)]$.  Following \cite{Cox}, we
regard $S_X$ as the semigroup algebra of $\NN^{\Sigma_{X}(1)}$ with
the $\Cl(X)$\nobreakdash-grading induced by $\deg(x^u) = \deg\bigl(
\prod_{\rho \in \Sigma_{X}(1)} x_{\rho}^{u_{\rho}} \bigr) = [u] \in
\Cl(X)$.  A divisor $D = \sum_{\rho \in \Sigma_{X}(1)}
u_{\rho}D_{\rho}$ determines a Laurent monomial $x^u = \prod_{\rho \in
  \Sigma_{X}(1)} x_{\rho}^{u_{\rho}} \in \kk[x_{\rho}^{\pm 1} : \rho
\in \Sigma_{X}(1)]$ and we often write the monomial as $x^D$.  The
support of $D$ or $x^u$ is the set
\[
\supp(D) = \supp(x^u) = \{ \rho \in \Sigma_{X}(1) : u_{\rho} \neq 0 \}
\, .
\]  
For a cone $\sigma \in \Sigma_{X}$, $\widehat{\sigma}$ is the set of
one-dimensional cones in $\Sigma_{X}$ that are not contained in
$\sigma$ and $x^{\widehat{\sigma}} = \prod_{\rho \in \widehat{\sigma}}
x_{\rho}$ is the associated monomial in $S_X$.  The irrelevant ideal
of $X$ is the square-free (i.e. reduced) monomial ideal $B_{X} := (
x^{\widehat{\sigma}} : \sigma \in \Sigma_{X})$.  Theorem~2.1 in
\cite{Cox} shows that the pair $(S_X,B_X)$ encodes a quotient
construction of $X$.  Specifically, if $\VV(B_X)$ is the subvariety of
$\AA^{\Sigma_X(1)}$ defined by $B_X$, then the toric variety $X$ is a
categorical quotient of $\AA^{\Sigma_{X}(1)} \setminus \VV(B_{X})$ by
the group $\Hom_{\ZZ}(\Cl(X), \kk^{*})$; the group action is induced
by the $\Cl(X)$-grading of $S_X$.

\subsection{Quivers}
\label{sub:quivers}

A quiver $Q$ is specified by two finite sets $Q_0$ and $Q_1$, whose
elements are called vertices and arrows, together with two maps
$\head, \tail \colon Q_1 \rTo Q_0$ indicating the vertices at the head
and tail of each arrow.  A nontrivial path in $Q$ is a sequence of
arrows $p = a_1 \dotsb a_\ell$ with $\head(a_{k}) = \tail(a_{k+1})$
for $1 \leq k < \ell$.  We set $\tail(p) = \tail(a_{1})$ and $\head(p)
= \head(a_\ell)$.  Each $i \in Q_0$ gives a trivial path $e_i$ where
$\tail(e_i) = \head(e_i) = i$.  The path algebra $\kk Q$ is the
$\kk$-algebra whose underlying $\kk$-vector space has a basis
consisting of paths in $Q$; the product of two basis elements equals
the basis element defined by concatenation of the paths if possible or
zero otherwise.  A cycle is a path $p$ in which $\tail(p) = \head(p)$.
A quiver is acyclic if it contains no cycles.  A vertex is a source if
it is not the head of any arrow and a quiver is rooted if it has a
unique source.

A walk $\gamma$ in $Q$ is an alternating sequence $i_0a_1i_1 \dotsb
a_\ell i_{\ell}$ of vertices $i_1, \dotsc, i_\ell$ and arrows
$a_1,\dotsc, a_\ell$ where $a_{k}$ is an arrow between $i_{k-1}$ and
$i_k$.  If $\tail(a_k) = i_{k-1}$ and $\head(a_k) = i_k$ then $a_k$ is
a forward arrow in $\gamma$; otherwise $\tail(a_k) = i_k$, $\head(a_k)
= i_{k-1}$ and $a_k$ is a backward arrow.  If $a \in Q_1$ then
$a^{-1}$ denotes the walk from $\head(a)$ to $\tail(a)$.  A walk
$\gamma$ is closed if $i_0 = i_\ell$ and a circuit is a closed walk in
which the arrows $a_1, \dotsc, a_\ell$ are distinct.  A quiver is
connected if there is a walk between any two vertices.  A tree is a
connected acyclic quiver.  We say that $Q' \subseteq Q$ is a spanning
subquiver if $Q_0' = Q_0$.

The vertex space $\ZZ^{Q_0}$ is the free abelian group of functions
from $Q_0$ to $\ZZ$ and the arrow space $\ZZ^{Q_1}$ is the free
abelian group of functions from $Q_1$ to $\ZZ$.  The characteristic
functions $\chi_{i} \colon Q_0 \rTo \ZZ$ for $i \in Q_0$ and $\chi_{a}
\colon Q_1 \rTo \ZZ$ for $a \in Q_1$ form the standard bases for the
vertex and arrows spaces.  We write $\NN^{Q_0}$ and $\NN^{Q_1}$ for
the semigroups generated by all $\NN$-linear combinations of the
characteristic functions $\chi_{i}$ and $\chi_{a}$ respectively.  The
incidence map $\inc \colon \ZZ^{Q_1} \rTo \ZZ^{Q_0}$ is defined by
$\inc(\chi_{a}) = \chi_{\head(a)} - \chi_{\tail(a)}$.  A function
$\theta \colon Q_0 \rTo \ZZ$ is an integral weight of $Q$ if $\sum_{i
  \in Q_0} \theta_i = 0$ and a function $f \colon Q_1 \rTo \ZZ$ is an
integral circulation if
\[
\sum\limits_{\substack{a \in Q_1\\ \tail(a) = i}} f_a =
\sum\limits_{\substack{a \in Q_1 \\ \head(a) = i}} f_a \;\; \text{for
  each $i \in Q_0$.}
\]
The weight lattice $\Wt(Q) \subset \ZZ^{Q_0}$ and the circulation
lattice $\Cir(Q) \subset \ZZ^{Q_1}$ are generated by the integral
weights and circulations respectively.  There is an exact sequence
\begin{equation} 
  \label{eq:graphses}
  \begin{diagram}[size=2.5em,l>=2.5em,midshaft,labelstyle=\scriptstyle]
    0 & \rTo & \Cir(Q) & \rTo & \ZZ^{Q_1} & \rTo^{\inc} & \Wt(Q) 
  \end{diagram}
\end{equation}
and the incidence map is surjective when $Q$ is connected; see \S4 in
\cite{Biggs}.  For $a \in Q_1$ and a walk $\gamma$, let
$\mult_{\gamma}(a) \in \ZZ$ equal the number of times $a$ appears as a
forward arrow in $\gamma$ minus the number of times it appears as a
backward arrow.  Given a walk $\gamma$, we set $f(\gamma) := \sum_{a
  \in Q_1} \mult_{\gamma}(a) \chi_{a} \in \ZZ^{Q_1}$; $f(\gamma) \in
\Cir(Q)$ if and only if $\gamma$ is a closed walk.

\subsection{Representations of Quivers} 
\label{sub:reps}

Let $Q$ be a connected quiver.  A representation $W = (W_i, w_a)$ of
$Q$ consists of a $\kk$-vector space $W_i$ for each $i \in Q_0$ and a
$\kk$-linear map $w_a \colon W_{\tail(a)} \rTo W_{\head(a)}$ for each
$a \in Q_1$. The dimension vector of $W$ is $\sum_{i\in Q_0}
\dim_{\kk}(W_{i})\chi_i \in \NN^{Q_0}$. In this paper, we will assume
that $\dim_{\kk}( W_{i}) = 1$ for all $i \in Q_0$.  A map between
representations $W = (W_i, w_a)$ and $W' = (W_i', w_a')$ is a family
$\psi_{i} \colon W_i^{\,} \rTo W_i'$ for $i \in Q_0$ of $\kk$-linear
maps that are compatible with the structure maps, that is $w_a'
\psi_{\tail(a)} = \psi_{\head(a)} w_a$ for all $a \in Q_1$.  With
composition defined componentwise, we obtain the abelian category of
representations of $Q$. Each rational weight $\theta \in
\Wt(Q)\otimes_{\ZZ}\QQ$ defines a stability notion for representations
and subquivers of $Q$.  A representation $W$ is $\theta$-stable if,
for every proper, nonzero subrepresentation $W' \subset W$, we have
$\theta(W') := \sum_{i \in \supp(W')} \theta_i > 0$, where $\supp(W')
:= \{ i \in Q_0 : W_i' \neq 0 \}$.  The notion of
$\theta$-semistability is obtained by replacing $>$ with $\geq$.  A
subquiver $Q' \subseteq Q$ is $\theta$-stable if it admits a
$\theta$-stable representation.

The isomorphism classes of representations are orbits in the
representation space
\[
\AA^{Q_1} = \Spec \bigl( \kk[y_a : a\in Q_1] \bigr) \cong
\bigoplus\limits_{a \in Q_1} \Hom_{\kk}(W_{\tail(a)}, W_{\head(a)})
\]
by the action of the group $(\kk^*)^{Q_0} \cong \prod_{i \in Q_0}
\operatorname{GL}(W_i)$ induced by the incidence map; in other words,
it acts by $(g \cdot w)_{a} = g_{\head(a)}^{\,} w_{a}
g_{\tail(a)}^{-1}$. Hence, we have a faithful action of the algebraic
torus $G := \Hom_{\ZZ} \bigl( \Wt(Q),\kk^* \bigr)$ on $\AA^{Q_1}$ and
it gives a $\Wt(Q)$-grading of the polynomial ring $S_Y := \kk[y_a :
a\in Q_1]$.  For $\theta \in \Wt(Q)$, let $(S_Y)_{\theta}$ be the
$\theta$-graded piece.  Following \cite{King1}, the GIT-quotient
\[
\mathcal{M}_\theta(Q) := \AA^{Q_1}\git G = \Proj \Biggl(
\bigoplus\limits_{k \in \NN} (S_Y)_{k\theta} \Biggr)
\]
is the categorical quotient $(\AA^{Q_1})^{\text{ss}}_\theta / G$,
where $(\AA^{Q_1})^{\text{ss}}_\theta \subseteq \AA^{Q_1}$ is the open
subscheme parametrizing $\theta$-semistable representations of $Q$.
Since $\dim(W_i) = 1$ for all $i \in Q_0$, $\mathcal{M}_\theta(Q)$ is
also a toric quiver variety as defined in \cite{Hille1}.  A weight
$\theta\in \Wt(Q)\otimes_{\ZZ}\QQ$ is generic if every
$\theta$-semistable representation is $\theta$-stable.  In this case,
$\mathcal{M}_\theta(Q)$ is the geometric quotient
$(\AA^{Q_1})^{\text{s}}_\theta / G$, where
$(\AA^{Q_1})^{\text{s}}_\theta$ parametrizes $\theta$-stable
representations of $Q$.  The set of generic weights decomposes into
finitely many open chambers, where $\mathcal{M}_\theta(Q)$ is
unchanged as $\theta$ varies in a chamber; see \cites{DolgachevHu,
  Thaddeus2}.

For generic $\theta$, Proposition~5.3 in \cite{King1} implies that
$\mathcal{M}_\theta(Q)$ is the fine moduli space of
$\theta$\nobreakdash-stable representations of $Q$.  To describe the
universal family on $\mathcal{M}_\theta(Q)$, we set $Q_0 = \{ 0,
\dotsc, r \}$ and identify the group $G$ with $\{ (g_0, \dotsc, g_r)
\in ( \kk^* )^{Q_0} : g_0 = 1 \}$.  This choice determines a
$G$-equivariant vector bundle $\bigoplus_{i \in Q_0}
\mathscr{O}_{\AA^{Q_1}}$ which descends to $\bigoplus_{i \in Q_0} F_i$
on $\mathcal{M}_\theta(Q)$; see Proposition~5.3 in \cite{King1}.  The
line bundles $F_0, \dotsc, F_r$ are called the tautological line
bundles on $\mathcal{M}_\theta(Q)$.  Since $G$ acts trivially on the
$0$-th component of $\bigoplus_{i \in Q_0} \mathscr{O}_{\AA^{Q_1}}$,
it follows that $F_0$ is the trivial line bundle.

\section{Quivers of Sections} 
\label{sec:quiver}

The goal of this section is to extend the classical notion of a linear
series from a single line bundle to a list of line bundles.  Let
$\mathcal{L} := (L_0, \dotsc, L_r)$ be a list of distinct line bundles
on the projective toric variety $X$.  A $T_X$-invariant section $s \in
H^0(X,L_j^{\;} \otimes L_{i}^{-1})$ is \emph{indecomposable} if the
divisor $\div(s)$ cannot be expressed as a sum $\div(s') + \div(s'')$
where $s' \in H^0(X,L_k^{\;} \otimes L_{i}^{-1})$ and $s'' \in
H^0(X,L_j^{\;} \otimes L_{k}^{-1})$ are nonzero $T_X$-invariant
sections and $0 \leq k \leq r$.  A \emph{quiver of sections}
associated to $\mathcal{L}$ is a quiver $Q$ in which the vertices $Q_0
= \{ 0, \dotsc, r \}$ correspond to the line bundles in $\mathcal{L}$
and the arrows from $i$ to $j$ correspond to a subset of the
indecomposable $T_X$-invariant sections in $H^0(X,L_j^{\;} \otimes
L_{i}^{-1})$.  If every indecomposable $T_X$-invariant section in
$H^0(X,L_j^{\;} \otimes L_{i}^{-1})$ for $0 \leq i, j \leq r$
corresponds to an arrow then $Q$ is the \emph{complete quiver of
  sections} for $\mathcal{L}$.  Since $X$ is projective, the unique
element in $H^0(X,L_i^{\;} \otimes L_{i}^{-1}) = H^0(X,
\mathscr{O}_{X})$ defines the trivial path $e_i$ in $Q$.  Moreover, if
$L_j \neq L_i$ then projectivity implies that both $H^{0}(X,L_j^{\;}
\otimes L_i^{-1})$ and $H^{0}(X,L_j^{-1} \otimes L_i^{\;})$ cannot be
nonzero.  It follows that $Q$ is acyclic.

\begin{conventions} 
  Let $Q$ be a quiver of sections associated to $\mathcal{L} = (L_0,
  \dotsc, L_r)$.
  \begin{enumerate}
  \item By definition, $Q$ only depends on the line bundles $L_j^{\,}
    \otimes L_{i}^{-1}$ where $0 \leq i, j \leq r$.  Consequently, for
    any line bundle $L'$ on $X$, we have $Q = Q'$ where $Q'$ is a
    quiver of sections associated to $\mathcal{L}' = (L_0 \otimes L',
    \dotsc, L_r \otimes L')$.  To eliminate this redundancy, we will
    assume that $L_0 = \mathscr{O}_{X}$.  By reordering the elements
    in $\mathcal{L}$ if necessary, we may also assume that $j < i$
    implies $H^0(X,L_j^{\;} \otimes L_{i}^{-1}) = 0$.
  \item We will assume that $H^0(X,L_i) \neq 0$ for $0 \leq i \leq r$.
    If $Q$ is the complete quiver of sections for $\mathcal{L}$, then
    this implies that $Q$ is connected and rooted at $0 \in Q_0$.
  \end{enumerate}  
\end{conventions}

Since each arrow $a \in Q_1$ corresponds to a $T_X$-invariant section
$s \in H^0(X,L_j^{\;} \otimes L_{i}^{-1})$, we simply write $\div(a)
:= \div(s) \in \CDiv(X)$.  More generally, for a path $p = a_1 \dotsb
a_\ell$ in $Q$, we set $\div(p) := \div(a_1) + \dotsb + \div(a_\ell)$.
This labelling of paths induces relations on $Q$.  Specifically, the
\emph{ideal of relations} is the two-sided ideal $R$ in the path
algebra $\kk Q$ generated by differences $p-p' \in \kk Q$ such that
$\tail(p) = \tail(p')$, $\head(p) = \head(p')$ and $\div(p) =
\div(p')$.  Since the arrows in $Q$ correspond to indecomposable
sections, $R$ is an admissible ideal.  The pair $(Q,R)$ is called a
\emph{bound quiver of sections};  the phrase ``bound quiver'' is a
synonym for ``quiver with relations''.

\begin{example}
  \label{exa:2pt}
  If $L_1$ is a nontrivial line bundle on $X$, then the complete
  quiver of sections for $\mathcal{L} = (\mathscr{O}_{X}, L_1)$ has two
  vertices and $\dim_{\kk} H^0(X,L_1)$ arrows.  The ideal of relations
  $R$ is the zero ideal.
\end{example}

The general correspondence between bound quivers and
finite-dimensional $\kk$-algebras has the following useful incarnation
for a quiver of sections.

\begin{proposition}
  \label{pro:algebra}
  If $(Q,R)$ is the complete bound quiver of sections for $\mathcal{L}
  \!=\! (L_0, \dotsc, L_r)$ then the quotient algebra $\kk Q/R$ is
  isomorphic to $\End\bigl( \bigoplus_{i=0}^{r} L_i \bigr)$.
\end{proposition}

\begin{proof}
  The map sending a path $p = a_1 \dotsb a_\ell$ in $Q$ to the product
  of the corresponding sections $s_1 \dotsb s_\ell \in
  H^0(X,L_{\head(p)}^{\;} \otimes L_{\tail(p)}^{-1})=
  \Hom(L_{\tail(p)}, L_{\head(p)})$ determines a homomorphism of
  $\kk$-algebras $\eta \colon \kk Q \rTo \End\bigl( \bigoplus_{i=0}^{r}
  L_i \bigr) = \bigoplus_{i,j=0}^{r} \Hom(L_i,L_j)$.  The map is
  surjective because $Q$ is a complete quiver.  Moreover, $\eta$ sends
  paths $p, p'$ in $Q$ satisfying $\tail(p) = \tail(p')$ and $\head(p)
  = \head(p')$ to the same element in $\Hom(L_{\tail(p)},
  L_{\head(p)})$ if and only if $\div(p) = \div(p')$.  Thus, we have
  $\Ker(\eta) = R$.
\end{proof}

\begin{remark}
  Let $(Q,R)$ be a (not necessarily complete) bound quiver of sections
  associated to $\mathcal{L} = (L_0, \dotsc, L_r)$.  If $\{ e_0,
  \dotsc, e_r \}$ is a complete set of primitive orthogonal
  idempotents in $\End\bigl( \bigoplus_{i=0}^{r} L_i \bigr)$ and $s_0,
  \dotsc, s_m$ are the indecomposable sections corresponding to the
  arrows in $Q$, then $\kk Q / R$ is the subalgebra of $\End\bigl(
  \bigoplus_{i=0}^{r} L_i \bigr)$ generated by $\{ e_0, \dotsc, e_r,
  s_0, \dotsc, s_m \}$.
\end{remark}

\begin{remark}
  \cites{Bondal, King1, BergmanProudfoot} work with the opposite
  quiver.  In particular, a ``Bondal quiver'' is a complete quiver of
  sections in which arrows have the opposite orientation.
\end{remark}

A quiver of sections $Q$ comes equipped with a distinguished lattice.
The map sending $a \in Q_1$ to $\div(a) \in \CDiv(X)$ extends to give
a $\ZZ$-linear map $\div \colon \ZZ^{Q_1} \rTo \CDiv(X)$ where $\div(v)
:= \sum_{a \in Q_1} v_{a} \div(a)$ for $v = \sum_{a \in Q_1} v_a
\chi_{a}$.  The \emph{section lattice} $\ZZ(Q)$ is the image of the
map $\pi := (\inc,\div) \colon \ZZ^{Q_{1}} \rTo \Wt(Q) \oplus
\CDiv(X)$; by definition, we have $\pi(\chi_a) = \bigl(
\chi_{\head(a)} - \chi_{\tail(a)}, \div(a) \bigr)$.  The projections
onto the components are denoted by $\pi_{1} \colon \ZZ(Q) \rTo \Wt(Q)$
and $\pi_{2} \colon \ZZ(Q) \rTo \CDiv(X)$ respectively.  These maps fit
in to the commutative diagram
\begin{equation}
  \label{eqn:rightsquare}
  \begin{diagram}[size=2em,midshaft,PS,labelstyle=\scriptstyle]
    \ZZ^{Q_1} & & & & \\
    & \rdTo\rdTo(4,2)^{\inc}\rdTo(2,4)_{\div} & & & \\
    & & \ZZ(Q) & \rTo_{\pi_1} & \Wt(Q) \\
    & & \dTo_{\pi_2} & & \dTo \\
    & & \CDiv(X) & \rTo & \Pic(X) \\
  \end{diagram}
\end{equation}
where $\pic(\theta) := \bigotimes_{i \in Q_0} L_{i}^{\theta_{i}}$ for
$\theta = \sum_{i \in Q_{0}} \theta_i \chi_{i}$.

\begin{example}
  \label{exa:FF1}
  Let $X = \FF_{1} = \PP\bigl( \mathscr{O}_{\PP^1} \oplus
  \mathscr{O}_{\PP^1}(1) \bigr)$ be the Hirzebruch surface determined
  by the fan in Figure~\ref{fig:FF1}~(a).
  \begin{figure}[!ht]
    \centering
    \mbox{
      \subfigure[Fan]{
        \psset{unit=1cm}
        \begin{pspicture}(-0.5,-0.5)(2.5,3.2)
          \pspolygon[linecolor=white,fillstyle=hlines*,
          hatchangle=25](1,1)(1,2)(2,2)(2,1)
          \pspolygon[linecolor=white,fillstyle=hlines*,
          hatchangle=5](1,1)(1,2)(0,2)(0,2)
          \pspolygon[linecolor=white,fillstyle=hlines*,
          hatchangle=55](1,1)(0,2)(0,0)(1,0)
          \pspolygon[linecolor=white,fillstyle=hlines*,
          hatchangle=165](1,1)(2,1)(2,0)(1,0)
          \cnode*[fillcolor=black](0,0){3pt}{P1}
          \cnode*[fillcolor=black](1,0){3pt}{P2}
          \cnode*[fillcolor=black](2,0){3pt}{P3}
          \cnode*[fillcolor=black](0,1){3pt}{P4}
          \cnode*[fillcolor=black](1,1){3pt}{P5}
          \cnode*[fillcolor=black](2,1){3pt}{P6}
          \cnode*[fillcolor=black](0,2){3pt}{P7}
          \cnode*[fillcolor=black](1,2){3pt}{P8}
          \cnode*[fillcolor=black](2,2){3pt}{P9}
          \ncline[linewidth=2pt]{->}{P5}{P7}
          \ncline{-}{P5}{P11}
          \ncline[linewidth=2pt]{->}{P5}{P6}
          \ncline[linewidth=2pt]{->}{P5}{P2}
          \ncline[linewidth=2pt]{->}{P5}{P8}
          \rput(2,0.7){\psframebox*{$1$}}
          \rput(1.3,2){\psframebox*{$2$}}
          \rput(-0.3,2){\psframebox*{$3$}}
          \rput(0.7,0){\psframebox*{$4$}}
        \end{pspicture}}
      \qquad \qquad  
      \subfigure[Quiver of sections]{
        \psset{unit=1.3cm}
        \begin{pspicture}(-0.5,-0.5)(2.5,2)
          \cnodeput(0,0){A}{0}
          \cnodeput(2,0){B}{1} 
          \cnodeput(0,2){C}{2}
          \psset{nodesep=0pt}
          \ncline[offset=5pt]{->}{A}{B} \lput*{:U}{$x_1$}
          \ncline[offset=-5pt]{->}{A}{B} \lput*{:U}{$x_3$}
          \ncline{->}{A}{C} \lput*{:270}{$x_4$}
          \ncline{<-}{C}{B} \lput*{:U}{$x_2$}
        \end{pspicture}}
      \qquad \qquad  
      \subfigure[Listing the arrows]{
        \psset{unit=1.3cm}
        \begin{pspicture}(-0.5,-0.5)(2.5,2)
          \cnodeput(0,0){A}{0}
          \cnodeput(2,0){B}{1} 
          \cnodeput(0,2){C}{2}
          \psset{nodesep=0pt}
          \ncline[offset=5pt]{->}{A}{B} \lput*{:U}{$a_1$}
          \ncline[offset=-5pt]{->}{A}{B} \lput*{:U}{$a_3$}
          \ncline{->}{A}{C} \lput*{:270}{$a_4$}
          \ncline{<-}{C}{B} \lput*{:U}{$a_2$}
        \end{pspicture}}
    }
    \caption{Hirzebruch surface $\FF_{1}$ \label{fig:FF1}}
  \end{figure}
  For $(k,\ell) \in \ZZ^2$, set $\mathscr{O}_{X}(k,\ell) :=
  \mathscr{O}_{X}(kD_1+\ell D_4) \in \Pic(X)$.  The complete quiver of
  sections for $\mathcal{L} = \bigl( \mathscr{O}_{X},
  \mathscr{O}_{X}(1,0), \mathscr{O}_{X}(0,1) \bigr)$ appears in
  Figure~\ref{fig:FF1}~(b).  The ideal of relations is $R = (0)$.  The
  section lattice $\ZZ(Q)$ is generated by the columns of the matrix
  \renewcommand{\arraystretch}{0.8}
  \renewcommand{\arraycolsep}{3pt}
  \begin{equation}
    \label{eq:FF1matrix}
    \left[ \text{\footnotesize $\begin{array}{rrrr}
        -1 & 0 & -1 & -1 \\
        1 & -1 & 1 & 0 \\
        0 & 1 & 0 & 1 \\
        1 & 0 & 0 & 0 \\
        0 & 1 & 0 & 0 \\
        0 & 0 & 1 & 0 \\
        0 & 0 & 0 & 1
      \end{array}$} \right] \, .
  \end{equation}
\end{example}

\begin{example}
  \label{exa:FF2}
  Let $X = \FF_{2} = \PP\bigl( \mathscr{O}_{\PP^1} \oplus
  \mathscr{O}_{\PP^1}(2) \bigr)$ be the Hirzebruch surface determined
  by the fan in Figure~\ref{fig:FF2}~(a).
  \begin{figure}[!ht]
    \centering
    \mbox{
      \subfigure[Fan]{
        \psset{unit=1cm}
        \begin{pspicture}(-0.5,-0.5)(2.5,3.2)
          \pspolygon[linecolor=white,fillstyle=hlines*,
          hatchangle=25](1,1)(1,3)(2,3)(2,1)
          \pspolygon[linecolor=white,fillstyle=hlines*,
          hatchangle=5](1,1)(1,3)(0,3)
          \pspolygon[linecolor=white,fillstyle=hlines*,
          hatchangle=55](1,1)(0,3)(0,0)(1,0)
          \pspolygon[linecolor=white,fillstyle=hlines*,
          hatchangle=165](1,1)(2,1)(2,0)(1,0)
          \cnode*[fillcolor=black](0,0){3pt}{P1}
          \cnode*[fillcolor=black](1,0){3pt}{P2}
          \cnode*[fillcolor=black](2,0){3pt}{P3}
          \cnode*[fillcolor=black](0,1){3pt}{P4}
          \cnode*[fillcolor=black](1,1){3pt}{P5}
          \cnode*[fillcolor=black](2,1){3pt}{P6}
          \cnode*[fillcolor=black](0,2){3pt}{P7}
          \cnode*[fillcolor=black](1,2){3pt}{P8}
          \cnode*[fillcolor=black](2,2){3pt}{P9}
          \cnode*[fillcolor=black](0,3){3pt}{P10}
          \cnode*[fillcolor=black](1,3){3pt}{P11}
          \cnode*[fillcolor=black](2,3){3pt}{P12}
          \ncline[linewidth=2pt]{->}{P5}{P10}
          \ncline{-}{P5}{P11}
          \ncline[linewidth=2pt]{->}{P5}{P6}
          \ncline[linewidth=2pt]{->}{P5}{P2}
          \ncline[linewidth=2pt]{->}{P5}{P8}
          \rput(2,0.7){\psframebox*{$1$}}
          \rput(1.3,2){\psframebox*{$2$}}
          \rput(-0.3,3){\psframebox*{$3$}}
          \rput(0.7,0){\psframebox*{$4$}}
        \end{pspicture}}
      \qquad \qquad  
      \subfigure[Quiver of sections]{
        \psset{unit=1.3cm}
        \begin{pspicture}(-0.5,-0.5)(2.5,2)
          \cnodeput(0,0){A}{0}
          \cnodeput(2,0){B}{1} 
          \cnodeput(0,2){C}{2}
          \cnodeput(2,2){D}{3}
          \psset{nodesep=0pt}
          \ncline[offset=5pt]{->}{A}{B} \lput*{:U}{$x_1$}
          \ncline[offset=-5pt]{->}{A}{B} \lput*{:U}{$x_3$}
          \ncline{->}{A}{C} \lput*{:270}{$x_4$}
          \ncline{->}{B}{D} \lput*{:270}{$x_4$}
          \ncarc[offset=1pt]{<-}{C}{B} \lput*{:U}{$x_1x_2$}
          \ncarc[offset=1pt]{->}{B}{C} \lput*{:180}{$x_2x_3$}
          \ncline[offset=5pt]{->}{C}{D} \lput*{:U}{$x_1$}
          \ncline[offset=-5pt]{->}{C}{D} \lput*{:U}{$x_3$}
        \end{pspicture}}\qquad \qquad  
      \subfigure[Listing the arrows]{
        \psset{unit=1.3cm}
        \begin{pspicture}(-0.5,-0.5)(2.5,2)
          \cnodeput(0,0){A}{0}
          \cnodeput(2,0){B}{1} 
          \cnodeput(0,2){C}{2}
          \cnodeput(2,2){D}{3}
          \psset{nodesep=0pt}
          \ncline[offset=5pt]{->}{A}{B} \lput*{:U}{$a_1$}
          \ncline[offset=-5pt]{->}{A}{B} \lput*{:U}{$a_2$}
          \ncline{->}{A}{C} \lput*{:270}{$a_3$}
          \ncline{->}{B}{D} \lput*{:270}{$a_6$}
          \ncarc[offset=1pt]{<-}{C}{B} \lput*{:U}{$a_4$}
          \ncarc[offset=1pt]{->}{B}{C} \lput*{:180}{$a_5$}
          \ncline[offset=5pt]{->}{C}{D} \lput*{:U}{$a_7$}
          \ncline[offset=-5pt]{->}{C}{D} \lput*{:U}{$a_8$}
        \end{pspicture}}
    }
    \caption{Hirzebruch surface $\FF_{2}$ \label{fig:FF2}}
  \end{figure}
  For $(k,\ell) \in \ZZ^2$, we write $\mathscr{O}_{X}(k,\ell) :=
  \mathscr{O}_{X}(kD_1+\ell D_4) \in \Pic(X)$.  The complete quiver of
  sections for $\mathcal{L} = \bigl( \mathscr{O}_X,
  \mathscr{O}_{X}(1,0), \mathscr{O}_{X}(0,1), \mathscr{O}_{X}(1,1)
  \bigr)$ appears in Figure~\ref{fig:FF2}~(b).  If we order the arrows
  as in Figure~\ref{fig:FF2}~(c), then the ideal of relations is $R =
  (a_2a_4-a_1a_5, a_4a_8-a_5a_7, a_2a_6-a_3a_8, a_1a_6-a_3a_7)$.  The
  section lattice $\ZZ(Q)$ is generated by the columns of the
  following matrix
  \renewcommand{\arraystretch}{0.8}
  \renewcommand{\arraycolsep}{2pt}
  \begin{equation}
    \label{eq:FF2matrix}
    \left[ \text{\footnotesize $\begin{array}{rrrrrrrr}
        -1 & -1 & -1 & 0 & 0 & 0 & 0 & 0 \\
        1 & 1 & 0 & -1 & -1 & -1 & 0 & 0 \\
        0 & 0 & 1 & 1 & 1 & 0 & -1 & -1 \\
        0 & 0 & 0 & 0 & 0 & 1 & 1 & 1 \\ 
        1 & 0 & 0 & 1 & 0 & 0 & 1 & 0 \\
        0 & 0 & 0 & 1 & 1 & 0 & 0 & 0 \\
        0 & 1 & 0 & 0 & 1 & 0 & 0 & 1 \\
        0 & 0 & 1 & 0 & 0 & 1 & 0 & 0
      \end{array}$} \right] \, ;
  \end{equation}
  the $i$-th column corresponds to $a_i$ for $1 \leq i \leq 8$.
\end{example}
 
Classically, one associates a projective space, called the linear
series, to a nonzero subspace of global sections of a line bundle.
Generalizing this construction, we associate a toric variety $Y_Q$ to
an appropriate quiver of sections $Q$.  Any connected, rooted, acyclic
quiver $Q$ defines a complete fan $\Sigma_Q$ in the $\RR$-vector space
$\Hom_{\ZZ}\bigl( \Cir(Q), \RR \bigr)$.  The rays $\Sigma_Q(1)$
correspond to arrows $a \in Q_1$ and are generated by the evaluation
maps $\eval_{a} \colon \Cir(Q) \rTo \RR$ defined by $\eval_a(f) = f_a$
for $f \in \Cir(Q)$.  Thus, $\eval_{a}$ is the unique generator of
$\rho_a \cap N_Y$ where $\rho_a \in \Sigma_Q(1)$ is the ray
corresponding to $a \in Q_1$ and $N_Y := \Hom_{\ZZ}\bigl( \Cir(Q), \ZZ
\bigr)$.  The rays $\rho_{a_1}, \dotsc, \rho_{a_\ell} \in \Sigma_Q(1)$
span a cone in $\Sigma_Q$ if and only if there exists a spanning tree
rooted at the source of $Q$ that does not contain $a_{1},\dots,
a_{\ell}$.  Hence, maximal cones in $\Sigma_Q$ correspond to spanning
trees rooted at the source, and have dimension $| Q_1|-| Q_0| + 1$.
Since $\Sigma_Q$ can also be described as the triangulation associated
to the region in the chamber complex corresponding to our acyclic
orientation of the underlying graph of $Q$ (combine Theorem~3.1 in
\cite{BGS} and Lemma~7.1 in \cite{GZ}), it follows that $\Sigma_Q$ is
a fan.  Let $Y_Q$ be the toric variety determined by the fan
$\Sigma_Q$.

The toric variety $Y_Q$, which is a toric quiver variety as defined in
\cite{Hille1}, has several other characterizations.  Following
\cite{BPS}, a toric variety $Y$ is unimodular if $M_Y$ is a unimodular
sublattice of $\ZZ^{\Sigma_Y(1)}$; see \eqref{eq:fundiagram}.  This is
equivalent to saying that $Y$ is smooth and any other variety obtained
from $Y$ by toric flips and flops is also smooth.

\begin{proposition}
  \label{pro:smooth}
  Let $Q$ be a connected, rooted, acyclic quiver.  If $Q_0 = \{ 0,
  \dotsc, r \}$ where $0$ is the unique source then the following
  varieties coincide:
  \begin{enumerate}
  \item the toric variety $Y_Q$ defined by the fan $\Sigma_Q$;
  \item the geometric quotient of $\AA^{Q_1} \setminus
    \mathbb{V}(B_Y)$ by the group $G = \Hom_{\ZZ}(\Wt(Q), \kk^{*})$,
    where $\AA^{Q_1} = \Spec(S_Y)$, $S_Y := \kk[y_a : a \in Q_1]$ and
    \[
    B_Y := \Biggl( \prod\limits_{a \in Q'_1} y_a : \text{$Q'$ is a
      spanning tree of $Q$ rooted at $0$} \Biggl) = \bigcap\limits_{i =
      1}^{r} \bigl( y_a : \head(a) = i \bigr) \, ;
    \]
  \item the fine moduli space $\mathcal{M}_\theta(Q)$ of
    $\theta$-stable representations for any rational weight $\theta
    \in \Wt(Q) \otimes_{\ZZ}\QQ$ lying in the open GIT-chamber
    containing $\vartheta:=\sum_{i\in Q_0} (\chi_i-\chi_0)$.
  \end{enumerate}
  Moreover,  this variety is unimodular and projective.
\end{proposition}

\begin{proof}
  Let $W$ be a $\vartheta$-semistable representation of $Q$.  If $W'
  \subset W$ is a proper nonzero subrepresentation, then we have
  $\vartheta(W') = \sum_{i \in \supp(W')} \vartheta_i \geq 0$.  Since
  $\vartheta_i = 1$ for $i \neq 0$ and $\vartheta_0 = -r$, it follows
  that $W_0' = 0$ and $\vartheta(W') > 0$.  Therefore, $\vartheta$ is
  generic, the open chamber containing $\vartheta$ is well-defined,
  and results in \S4 of \cite{King1} show that
  $\mathcal{M}_\vartheta(Q)$ is a smooth projective variety.  Since
  the map $\inc\colon \ZZ^{Q_1}\rTo \Wt(Q)$ is totally unimodular
  (e.g. Proposition~5.3 in \cite{Biggs} or Example~2 in \S19.3 of
  \cite{Sch}), it follows that $\Cir(Q)$ is a unimodular sublattice of
  $\ZZ^{Q_1}$ and the toric variety defined by $\Sigma_Q$ is
  unimodular. Theorem~2.1 in \cite{Cox} establishes the equivalence
  between (a) and (b), and the discussion preceding Theorem~1.7 in
  \cite{Hille1} establishes the equivalence between (a) and (c).
\end{proof}

When $Q$ is a connected, rooted, acyclic quiver of sections on $X$,
the toric variety $Y_Q$ is called the \emph{multilinear series} of
$Q$.  When the quiver of sections is unambiguous, we simply write $Y$
for the multilinear series.  If $Q$ is the complete quiver of sections
for a list a of line bundles $\mathcal{L}$, then we write $|
\mathcal{L}| := Y$ for the \emph{complete multilinear series}.

\begin{remark}
  \label{rem:Y}
  Proposition~\ref{pro:smooth} implies that the diagram
  \eqref{eq:fundiagram} becomes
  \begin{equation*} 
    \begin{diagram}[size=2em,l>=2em,midshaft,labelstyle=\scriptstyle]
      0 & \rTo & M_Y & \rTo & \CDiv(Y) & \rTo & \Pic(Y) & \rTo & 0 \,
      \, \\
      & & \dTo & & \dTo & & \dTo & & \\
      0 & \rTo & \Cir(Q) & \rTo & \ZZ^{Q_1} & \rTo^{\inc} & \Wt(Q) &
      \rTo & 0 \, ,
    \end{diagram}
  \end{equation*}
  where the vertical maps are isomorphisms.  The open chamber
  containing $\vartheta$ is the ample cone $\Amp_\QQ(Y)$ of
  $\QQ$\nobreakdash-divisor classes on $Y$ and the closure of this
  chamber is
  \[
  \Nef_\QQ(Y) = \bigcap\limits_{Q' \subseteq Q} \Biggl\{ \sum\limits_{a
    \in Q_1'} \lambda_a [D_a] : \lambda_a \in \QQ_{\geq 0} \Biggr\} \,
  ,
  \]
  where the intersection runs over all spanning trees $Q'$ of $Q$
  rooted at $0$ and $D_a$ is the irreducible $T_Y$-invariant Weil
  divisor associated to $a \in Q_1$.  Since $Y$ is smooth, the ample
  line bundle $\mathscr{O}_{Y}(\vartheta)$ determined by $\vartheta
  \in \Wt(Q)$ is very ample.
\end{remark}

\begin{example}
  \label{exa:2ptY}
  Let $Q$ be a quiver with $Q_0 = \{ 0, 1 \}$ and $Q_1 = \{ a_0,
  \dotsc, a_m \}$ such that $\tail(a_k) = 0$ and $\head(a_k) = 1$ for
  all $0 \leq k \leq m$.  Since every arrow forms a spanning tree
  rooted at $0$, the irrelevant ideal of $Y$ is $B_Y = ( y_{a_0},
  \dotsc, y_{a_{m}})$.  Hence, $Y$ is the geometric quotient of
  $\AA^{Q_1} \setminus \{ 0 \}$ by $G := \Hom_{\ZZ}\bigl( \Wt(Q),
  \kk^{*})$.  Choosing $\chi_1 - \chi_0$ as a basis for $\Wt(Q)$, we
  see that $G \cong \kk^*$ and the $G$-action is induced by the matrix
  is $\left[ \begin{smallmatrix} 1 & \dotsb & 1 \end{smallmatrix}
  \right]$.  Therefore, we have $Y \cong \PP^{m}$.  In particular, if
  $Q$ is the complete quiver of sections for $\mathcal{L} =
  (\mathscr{O}_{X}, L_1)$ described in Example~\ref{exa:2pt}, then the
  complete multilinear series $|\mathcal{L}|$ is canonically
  isomorphic to the linear series $|L_1|$.
\end{example}

\begin{example}
  \renewcommand{\arraystretch}{0.8}
  \renewcommand{\arraycolsep}{3pt}
  \label{exa:FF1Y}
  Let $X = \FF_1$ and $\mathcal{L} = \bigl( \mathscr{O}_X,
  \mathscr{O}_{X}(1,0), \mathscr{O}_{X}(0,1) \bigr)$ as in
  Example~\ref{exa:FF1}.  If we identify $\Cir(Q)$ with $\ZZ^2$ by
  choosing the circuits $(a_1^{\,}a_3^{-1}, a_3^{\,}a_2^{\,}a_4^{-1})$
  as an ordered basis, then the unique generator of $\rho_{k} \cap
  N_Y$, where $\rho_k \in \Sigma_Q(1)$ corresponds to $a_k \in Q_1$,
  is the $k$-th column of the matrix $\bigl[
    \text{\footnotesize $\begin{array}{rrrr}
      1 & 0 & -1 & 0 \\ 0 & 1 & 1 & -1
    \end{array}$} \bigr]$.
  Figure~\ref{fig:FF1}~(c) gives $S_Y = \kk[y_1, \dotsc, y_4]$ and
  $B_Y = (y_1,y_3) \cap (y_2, y_4)$.  Hence, the quotient construction
  of $Y$ from Proposition~\ref{pro:smooth}~(b) coincides with the
  quotient construction of $X$ encoded by the pair $(S_X,B_X)$; see
  \S\ref{sub:toric}.  Therefore, the multilinear series $Y$ equals
  $X$.
\end{example}

\begin{example}
  \label{exa:FF2Y}
  Let $X = \FF_2$ and $\mathcal{L} = \bigl( \mathscr{O}_X,
  \mathscr{O}_{X}(1,0), \mathscr{O}_{X}(0,1), \mathscr{O}_{X}(1,1)
  \bigr)$ as in Example~\ref{exa:FF2}.  If we identify $\Cir(Q)$ with
  $\ZZ^5$ by choosing the circuits
  \[
  (a_1^{\,}a_2^{-1}, a_1^{\,}a_4^{\,}a_3^{-1},
  a_1^{\,}a_5^{\,}a_3^{-1}, a_1^{\,}a_6^{\,}a_7^{-1}a_3^{-1},
  a_1^{\,}a_6^{\,}a_8^{-1}a_3^{-1})
  \] 
  as an ordered basis, then the unique generator of $\rho_{k} \cap
  N_Y$, where $\rho_k \in \Sigma_Q(1)$ corresponds to $a_k \in Q_1$,
  is the $k$-th column of the matrix
  \renewcommand{\arraystretch}{0.8}
  \renewcommand{\arraycolsep}{3pt}
  \[
  \left[ 
    \text{\footnotesize $\begin{array}{rrrrrrrr}
      1 & -1 & 0 & 0 & 0 & 0 & 0 & 0 \\
      1 & 0 & -1 & 1 & 0 & 0 & 0 & 0 \\
      1 & 0 & -1 & 0 & 1 & 0 & 0 & 0 \\
      1 & 0 & -1 & 0 & 0 & 1 & -1 & 0 \\
      1 & 0 & -1 & 0 & 0 & 1 & 0 & -1 
    \end{array}$} \right] \, .
  \]
  Figure~\ref{fig:FF2}~(c) implies that $S_Y = \kk[y_1, \dotsc, y_8]$
  and $B_Y = (y_1,y_2) \cap (y_3, y_4, y_5) \cap (y_6, y_7, y_8)$.
  Hence, the multilinear series $Y$ is a smooth $5$-dimensional toric
  variety with $8$ irreducible $T_Y$-invariant Weil
  divisors and $18$ $T_Y$-fixed points.  The ample cone of $Y$ is 
  \[
  \Amp_\QQ(Y) = \left\{ \theta = (\theta_0, \theta_1,
    \theta_2, \theta_3) \in \Wt(Q) \otimes_{\ZZ} \QQ : \text{$\theta_1
      > 0$, $\theta_2 > 0$, $\theta_3 > 0$} \right\} \, .
  \]
  Since $(-3,-1,1,3) \not\in \Amp_\QQ(Y)$, the dualizing line bundle
  on $Y$ is not ample.
\end{example}

\section{Multilinear Series}
\label{sec:multilinear}

In this section, we study morphisms from the toric variety $X$ to the
multilinear series $Y$ induced by a quiver of sections.  To begin, we
give necessary and sufficient conditions for a quiver of sections $Q$
on $X$ to define a morphism from $X$ to the multilinear series $Y$.
The divisors labelling the arrows in $Q$ define a ring homomorphism
$\Phi_Q \colon S_{Y} \rTo S_X$ between the total coordinate rings of
$X$ and $Y$ given by $\Phi_Q(y_a) = x^{\div(a)}$.  The \emph{base
  ideal} of $Q$ is the ideal $B_Q$ in $S_X$ generated by the image
$\Phi_Q(B_Y)$.

\begin{proposition} 
  \label{pro:free} 
  Let $Q$ be a connected, rooted, acyclic quiver of sections on $X$.
  If $Y$ is the multilinear series of $Q$, then the following are
  equivalent:
  \begin{enumerate}
  \item the map $\Phi_Q$ determines a morphism $\varphi_Q
    \colon X \rTo Y$;
  \item the irrelevant ideal $B_X$ is contained in the radical
    $\rad(B_Q)$;
  \item for all $\sigma \in \Sigma_X$, there exists a spanning tree
    $Q' \subseteq Q$ rooted at the unique source in $Q$ such that
    $\supp(\div(a)) \subseteq \widehat{\sigma}$ for all $a \in Q'_1$.
  \end{enumerate}
\end{proposition}

\begin{proof}
  The toric variety $X$ is a categorical quotient of
  $\AA^{\Sigma_{X}(1)} \setminus \VV(B_X)$ under the action of the
  group $\Hom_{\ZZ}(\Cl(X), \kk^*)$; see \S\ref{sub:toric}.
  Similarly, Proposition~\ref{pro:smooth} shows that $Y$ is the
  geometric quotient of $\AA^{Q_1} \setminus \VV(B_Y)$ under the
  action of the group $G = \Hom_{\ZZ}(\Pic(Y), \kk^*)$.  The ring map
  $\Phi_Q \colon S_Y \rTo S_X$ defines a morphism from
  $\AA^{\Sigma_{X}(1)}$ to $\AA^{Q_1}$.  This morphism is equivariant
  with respect to the actions of the groups $\Hom_{\ZZ}(\Cl(X),
  \kk^*)$ and $G$ on $\AA^{\Sigma_{X}(1)}$ and $\AA^{Q_1}$
  respectively because the lattice maps
  \[
  \begin{diagram}[size=2em,l>=2em,midshaft,labelstyle=\scriptstyle]
    \ZZ^{Q_1} & \rTo^{\inc} & \Wt(Q) \\
    \dTo^{\div} & &  \dTo_{[\pic]} \\
    \ZZ^{\Sigma_X(1)} & \rTo & \Cl(X)  \,     
  \end{diagram}
  \]
  commute.  Thus, $\Phi_Q$ induces the morphism $\varphi_Q \colon X
  \rTo Y$ if and only if the preimage of the irrelevant set $\VV(B_Y)$
  is contained in the irrelevant set $\VV(B_X)$; see Theorem~3.2 in
  \cite{Cox2}.  Since the preimage of $\VV(B_Y)$ is cut out by the
  base ideal $B_Q$, this is equivalent to $B_X$ being contained in
  $\rad(B_Q)$.  In other words, conditions (a) and (b) are equivalent.
  The definition of $B_X$ and the explicit description of $B_Y$ from
  Proposition~\ref{pro:smooth}~(b) gives the equivalence between (b)
  and (c).
\end{proof}

A quiver of sections $Q$ is \emph{basepoint-free} if it is connected,
rooted, acyclic, and satisfies any of the equivalent conditions in
Proposition~\ref{pro:free}.  If the complete quiver of sections for
$\mathcal{L}$ is basepoint-free, then $\varphi_{|\mathcal{L}|} \colon
X \rTo |\mathcal{L}|$ denotes the associated morphism to the complete
multilinear series.

\begin{corollary}
  \label{cor:basepointfree}
  If $Q$ is a basepoint-free quiver of sections then each line bundle
  $L_i$ on $X$ is basepoint-free.  Conversely, if each $L_i$ is
  basepoint-free and $Q$ is a complete quiver of sections for
  $\mathcal{L} = (\mathscr{O}_X,L_1, \dotsc, L_r)$ then $Q$ is
  basepoint-free.
\end{corollary}

\begin{proof} 
  Since $Q$ is basepoint-free, it satisfies condition (c) from
  Proposition~\ref{pro:free}.  Hence, for $i \in Q_0$ and $\sigma \in
  \Sigma_X(d)$, there exists a path $p = a_1 \dotsc a_\ell$ in $Q$
  such that $\tail(p) = 0$, $\head(p) = i$, and
  $\supp\bigl(\div(a_k)\bigr) \subseteq \widehat{\sigma}$ for all $1
  \leq k \leq \ell$.  In other words, $L_i$ admits a
  $T_X$\nobreakdash-invariant section that does not vanish at the
  $T_X$-fixed point indexed by $\sigma$ for all $\sigma \in
  \Sigma_X(d)$.  Therefore, for each $i \in Q_0$, $L_i$ is
  basepoint-free.  When $Q$ is complete, we can reverse
  the argument for the first part.
\end{proof}

Given a basepoint-free quiver of sections $Q$, the image of
$\varphi_Q$ can be described explicitly.  Let $\NN(Q)$ be the image of
$\NN^{Q_1}$ under the map $\pi \colon \ZZ^{Q_1} \rTo \Wt(Q) \oplus
\CDiv(X)$.  Observe that $\NN(Q)$ is a subsemigroup of the section
lattice $\ZZ(Q)$.  Since $S_Y$ is the semigroup algebra of
$\NN^{Q_1}$, the map $\pi$ induces a surjective $\kk$-algebra
homomorphism from $S_Y$ to $\kk[\NN(Q)]$ where $\kk[\NN(Q)]$ is the
semigroup algebra of $\NN(Q)$.  The kernel of this induced map is the
toric ideal $I_Q := \bigl( y^u - y^v \in S_Y : u - v \in \Ker(\pi)
\bigr)$.  The $G$-invariant affine toric variety $\VV(I_Q) \subseteq
\AA^{Q_1}$ cut out by the ideal $I_Q$ need not be normal.  The ideal
$I_Q$ is analogous to the toric ideal defined by the augmented
vertex-edge incidence matrix of the McKay quiver in \cite{CMT2}.  With
this notation, we obtain the following.

\begin{proposition}
  \label{pro:ideal}
  Let $Q$ be a basepoint-free quiver of sections on $X$.  If
  $\varphi_Q \colon X \rTo Y$ is the induced morphism, then the image
  of $\varphi_Q$ is:
  \begin{enumerate}
  \item the subscheme of $Y$ corresponding to the $\Wt(Q)$-graded
    $B_Y$-saturated ideal $I_Q$;
  \item the geometric quotient of $\VV(I_Q) \setminus \VV(B_Y)$ by $G
    = \Hom_\ZZ\bigl( \Wt(Q), \kk^* \bigr)$;
  \item the GIT-quotient $\VV(I_Q) \git G$ for any $\theta \in
    \Amp_\QQ(Y)$.
  \end{enumerate}
\end{proposition}

\begin{proof}
  Proposition~\ref{pro:smooth}~(b) implies that the closed subsets of
  $Y$ are in bijection with the $G$-invariant closed subsets of
  $\AA^{Q_1} \setminus \VV(B_Y)$, and hence with the $B_Y$-saturated
  ideals of $S_Y$.  The image of the map from
  $\mathbb{A}^{\Sigma_X(1)}$ to $\mathbb{A}^{Q_1}$ induced by the
  homomorphism of semigroups $\div \colon \mathbb{N}^{Q_1} \rTo
  \CDiv(X)$ is cut out by the toric ideal $\Ker(\Phi_Q)$.  Since the
  action of $G = \Hom_{\ZZ}\bigl(\Wt(Q),\kk^*\bigr)$ on $\AA^{Q_1}$ is
  induced by the map $\inc \colon \ZZ^{Q_1} \rTo \Wt(Q)$, the toric
  ideal $I_Q$ associated to the map $\pi = (\inc,\div) \colon
  \NN^{Q_1} \rTo \Wt(Q) \oplus \CDiv(X)$ cuts out the
  $\Wt(Q)$\nobreakdash-homogeneous part of $\Ker(\Phi_Q)$.  As $I_Q$
  is prime and hence $B_Y$\nobreakdash-saturated, (a) and (b) follow.
  The equivalence of (b) and (c) follows directly from the equivalence
  between (b) and (c) in Proposition~\ref{pro:smooth}.
\end{proof}

\begin{remark}
  \label{rem:chamber}
  Proposition~\ref{pro:ideal}~(c) holds for rational weights in a cone
  that may strictly contain $\Amp_\QQ(Y)$ because the GIT-chamber
  decomposition for the $G$-action on $\VV(I_Q)$ is a coarsening of
  that for the $G$-action on $\AA^{Q_1}$.
\end{remark}

\begin{proof}[Proof of Theorem~\ref{thm:1}]
  The complete quiver of sections for $\mathcal{L} = (\mathscr{O}_X,
  L_1, \dotsc, L_r)$ is basepoint-free by
  Corollary~\ref{cor:basepointfree}, so $\varphi_{|\mathcal{L}|}$ is a
  morphism.  The explicit description of the image as a geometric
  quotient is presented in Proposition~\ref{pro:ideal}~(b).
\end{proof}

\begin{example}
  \label{exa:2ptimage}
  Corollary~\ref{cor:basepointfree} show that the complete quiver of
  sections for $\mathcal{L} = (\mathscr{O}_X, L_1)$ is basepoint-free
  if and only if $L_1$ is basepoint-free.  Since the semigroup
  $\NN(Q)$ is generated by $\bigl( (-1, 1), \div(s) \bigr) \in \ZZ^2
  \oplus \CDiv(X)$ for $T_X$-invariant sections $s \in H^0(X, L_1)$,
  it is isomorphic to the semigroup generated by the effective
  divisors $\div(s)$ in $\CDiv(X)$.  Hence $I_Q$ is the ideal for the
  image of $X$ in $|L_1| := \PP \bigl( H^0(X,L_1) \bigr)$.  It follows
  that $\VV(I_Q)$ is the affine cone over $\varphi_{|L_1|}(X)$, and
  $\varphi_{|\mathcal{L}|} = \varphi_{|L_1|}$.
\end{example}

\begin{example}
  \label{exa:FF1image}
  Let $X = \FF_1$ and $\mathcal{L} = \bigl( \mathscr{O}_X,
  \mathscr{O}_{X}(1,0), \mathscr{O}_{X}(0,1) \bigr)$ as in
  Example~\ref{exa:FF1}.  Example~\ref{exa:FF1Y} shows that the
  multilinear series is $Y = |\mathcal{L}| = X$.  By examining
  Figure~\ref{fig:FF1}~(b), we see that the complete quiver of
  sections for $\mathcal{L}$ satisfies condition~(c) in
  Proposition~\ref{pro:free}.  By definition, $\NN(Q)$ is generated by
  the columns of the matrix in \eqref{eq:FF1matrix} so the toric ideal
  is $I_Q = (0) \subset S_Y$.  Thus, $\varphi_{|\mathcal{L}|}$ is an
  isomorphism.
\end{example}

\begin{example}
  \label{exa:FF2image}
  Let $X = \FF_2$ and $\mathcal{L} = \bigl( \mathscr{O}_X,
  \mathscr{O}_{X}(1,0), \mathscr{O}_{X}(0,1), \mathscr{O}_{X}(1,1)
  \bigr)$ as in Example~\ref{exa:FF2}.  The multilinear series $Y =
  |\mathcal{L}|$ is described in Example~\ref{exa:FF2Y}.  By examining
  Figure~\ref{fig:FF2}~(b), we see that the complete quiver of
  sections for $\mathcal{L}$ satisfies condition~(c) in
  Proposition~\ref{pro:free}.  By definition, $\NN(Q)$ is generated by
  the columns of the matrix in \eqref{eq:FF2matrix}, so $I_Q =
  (y_2y_4-y_1y_5, y_1y_6-y_3y_7, y_2y_6-y_3y_8, y_2y_7-y_1y_8,
  y_5y_7-y_4y_8) \subset S_Y$.  In this case, we have
  $\varphi_{|\mathcal{L}|}(X) = \VV(I_Q) \git G$ for all rational
  weights $\theta$ in the GIT-chamber
  \[
  \Theta := \bigl\{ (\theta_0,\theta_1,\theta_2,\theta_3) \in \Wt(Q)
  \otimes_\ZZ \QQ : 
  \theta_3 > 0,\; \theta_1 + \theta_3 > 0,\; \theta_2 + \theta_3>0,\;
  \theta_1+\theta_2+\theta_3>0 \bigr\} \, ;
  \]
  see Remark~\ref{rem:chamber}.
\end{example}

With additional hypotheses, we can enlarge the commutative diagram
\eqref{eqn:rightsquare}.

\begin{corollary}
  \label{cor:funddiagram}
  If $Q$ is basepoint-free and $\dim \varphi_Q(X) = \dim X$, then
  \begin{equation*}
    \begin{diagram}[size=2em,l>=2em,midshaft,labelstyle=\scriptstyle]
      0 & \rTo & M_X & \rTo & \ZZ(Q) & \rTo^{\pi_1} & \Wt(Q) & \rTo &
      0 \\
      & & \dEqual & & \dTo_{\pi_2} & & \dTo_{\pic} & & \\
      0 & \rTo & M_X & \rTo & \CDiv(X) & \rTo & \Pic(X) & \rTo & 0 \\ 
    \end{diagram}
  \end{equation*} 
  is a commutative diagram with exact rows.  In particular, the
  projection $\pi_2$ induces an isomorphism between $\Ker(\pi_1)$ and
  $M_X$.
\end{corollary}

\begin{proof}
  Combining \eqref{eqn:rightsquare} with the top row of
  \eqref{eq:fundiagram}, it is enough to prove that $\pi_2$ yields an
  isomorphism between $\Ker(\pi_1)$ and $M_X$.  The morphism
  $\varphi_Q \colon X \rTo \varphi_Q(X)$ corresponds to the map of
  semigroup algebras $S_Y/I_Q = \kk[\mathbb{N}(Q)] \rTo
  \kk[\mathbb{N}^{\Sigma_X(1)}] = S_X$ induced by $\pi_2$.  Since
  $\dim \varphi_Q(X) = \dim X$, it identifies the dense tori in $X$
  and $\varphi_Q(X)$.  Therefore, $\pi_2$ identifies the character
  lattices $M_X$ and $\Ker(\pi_1)$.
\end{proof}

Next, we give a criterion for $\varphi_Q \colon X \rTo Y$ to be a
closed embedding.  For $\sigma \in \Sigma_X(d)$, let
$y^{\widehat{\sigma}} := \prod_{\supp(\div(a)) \subseteq
  \widehat{\sigma}} y_a$ be the associated monomial in $S_Y$.  The
localization of an $S_Y$-module $F$ at the element $y^{-
  \widehat{\sigma}}$ is denoted by $F[y^{- \widehat{\sigma}}]$.  The
weight $\vartheta:=\sum_{i\in Q_0} (\chi_i-\chi_0)$ appearing below is
defined in Proposition~\ref{pro:smooth}~(c).

\begin{proposition}
  \label{pro:ample}
  Let $Q$ be a basepoint-free quiver of sections.  The map $\varphi_Q
  \colon X \rTo Y$ is a closed embedding if and only if the line bundle
  $L := L_0^{\vartheta_0} \otimes \dotsb \otimes L_r^{\vartheta_r} =
  \bigotimes_{i \in Q_0} L_i$ is ample and $\bigl( (S_Y/I_Q)[y^{-
    \widehat{\sigma}}] \bigr)_{[0]} \cong \bigl( S_X[x^{-
    \widehat{\sigma}}] \bigr)_{[0]}$ for all $\sigma \in \Sigma_X(d)$.
\end{proposition}

\begin{proof}
  The very ample line bundle $\mathscr{O}_{Y}(\vartheta)$ from
  Remark~\ref{rem:Y} provides a closed embedding $Y \rTo \PP^{m} := \PP
  \bigl( H^0(Y, \mathscr{O}_{Y}(\vartheta)) \bigr)$.  Hence,
  $\varphi_Q \colon X \rTo Y$ is a closed embedding if and only if the
  composition $X \rTo \PP^m$ is a closed embedding.  The map $X \rTo
  \PP^m$ is determined by $\varphi_Q^{*}\bigl(
  \mathscr{O}_{Y}(\vartheta) \bigr) = \pic(\vartheta)= L$ and the
  subspace $\Phi_Q\bigl( (S_Y)_\vartheta \bigr) \subseteq (S_X)_{[L]}
  \cong H^0( X, L)$ of global sections.  The morphism $\varphi_Q
  \colon X \rTo Y$ corresponds to the map of semigroup algebras
  $S_Y/I_Q = \kk[\mathbb{N}(Q)] \rTo \kk[\mathbb{N}^{\Sigma_X(1)}] =
  S_X$ induced by $\pi_2 \colon \ZZ(Q) \rTo \CDiv(X)$ which implies
  that $(S_Y/I_Q)_\vartheta \cong \Phi_Q((S_Y)_\vartheta)$.  Moreover,
  the set $\mathcal{V} := \NN(Q)\cap \pi_1^{-1}(\vartheta)$ can be
  identified with both the monomial basis of $(S_Y/I_Q)_\vartheta$ and
  a subset of lattice points in the polytope associated to $L$.  Since
  $Q$ is basepoint-free, Proposition~\ref{pro:free}~(c) implies that,
  for each cone $\sigma \in \Sigma_X(d)$, there exists a monomial
  $y^{u_{\sigma}} \in (S_Y)_\vartheta$ such that $x^{v_\sigma} :=
  \Phi_Q(y^{u_\sigma}) \in (S_X)_{[L]}$ satisfies $\supp(x^{v_\sigma})
  \subseteq \widehat\sigma$.  By Theorem~2.7 in \cite{Oda}, the
  monomial $x^{v_\sigma} \in S_X$ corresponds to the vertex $v_\sigma
  \in \mathcal{V}$ of the polytope associated to $L$.  From the proof
  of the Theorem~2.13 in \cite{Oda}, we deduce that $X \rTo \PP^m$ is a
  closed embedding if and only if $v_\sigma \neq v_\tau$ holds for
  each pair $\sigma \neq \tau \in \Sigma_X(d)$, and the semigroup $M_X
  \cap \sigma^{\vee}$ is generated by $\mathcal{V} - v_\sigma$ for all
  $\sigma \in \Sigma_X(d)$.  Corollary~2.14 in \cite{Oda} proves that
  the first condition is equivalent to $L$ being ample.  Lemma~2.2 in
  \cite{Cox} shows that $\kk[M_X \cap \sigma^{\vee}] \cong
  (S_X[x^{-\widehat\sigma}])_{[0]}$.  Working in the semigroup algebra
  $\kk[\NN(Q)] = S_Y/I_Q$ rather than $\kk[\NN^{\Sigma_X(1)}] = S_X$,
  essentially the same argument establishes that $\bigl(
  (S_Y/I_Q)[y^{- \widehat\sigma}] \bigr)_{[0]}$ is isomorphic to the
  semigroup algebra of $\mathcal{V} - v_{\sigma}$.  Therefore, the
  second condition is equivalent to $\bigl( (S_Y/I_Q)[y^{-
    \widehat{\sigma}}] \bigr)_{[0]} \cong \bigl( S_X[x^{-
    \widehat{\sigma}}] \bigr)_{[0]}$ for all $\sigma \in \Sigma_X(d)$.
\end{proof}

A quiver of sections $Q$ is \emph{very ample} if it is complete,
basepoint-free, and $\varphi_Q \colon X \rTo Y$ is a closed embedding.
For convenience, we record an instance of Proposition~\ref{pro:ample}.

\begin{corollary}
  \label{cor:veryample}
  Let $\mathcal{L} = (\mathscr{O}_X, L_1, \dotsc, L_r)$ a list of
  basepoint-free line bundles and set $L := \bigotimes_{i \in Q_0}
  L_i$.  Assume that the map $H^0(X,L_1) \otimes_\kk \dotsb
  \otimes_\kk H^0(X,L_r) \rTo H^0(X,L)$ is surjective.  The morphism
  $\varphi_{|\mathcal{L}|} \colon X \rTo Y$ is a closed embedding if
  and only if $L$ is very ample.
\end{corollary}

\begin{proof}
  Since each $L_i$ is basepoint-free,
  Corollary~\ref{cor:basepointfree} implies that the complete quiver
  of sections $Q$ for $\mathcal{L}$ is basepoint-free, so
  $(S_Y/I_Q)_{\chi_i - \chi_0} \cong (S_X)_{[L_i]} \cong H^0(X,L_i)$.
  Since $L_0 = \mathscr{O}_X$ and $\vartheta = \sum_{i=1}^{r} (\chi_i
  - \chi_0)$, we may identify the image of the map
  \[
  H^0(X,L_0^{\vartheta_0}) \otimes_\kk \dotsb \otimes_\kk
  H^0(X,L_r^{\vartheta_r}) = H^0(X,L_1) \otimes_\kk \dotsb \otimes_\kk
  H^0(X,L_r) \rTo H^0(X,L)
  \]
  with the vector space $(S_Y/I_Q)_\vartheta$.  From the proof of
  Proposition~\ref{pro:ample}, we know that $\varphi_Q \colon X \rTo Y$
  is a closed embedding if and only if the map to projective space,
  determined by the line bundle $L$ and the subspace
  $(S_Y/I_Q)_\vartheta \subseteq H^0(X,L)$, is a closed embedding.
  The hypothesis that $H^0(X,L_1) \otimes_\kk \dotsb \otimes_\kk
  H^0(X,L_r) \rTo H^0(X,L)$ is surjective implies that
  $(S_Y/I_Q)_\vartheta \cong H^0(X,L)$.  Lastly, we observe that the
  complete linear series $|L|$ determines a closed embedding if and
  only $L$ is very ample.
\end{proof}

\begin{example}
  Corollary~\ref{cor:veryample} shows that the complete quiver of
  sections for the list $\mathcal{L} = (\mathscr{O}_X, L_1)$ is very
  ample if and only the line bundle $L_1$ is very ample.
\end{example}

\begin{example}
  \label{exa:FF1ample}
  Let $X = \FF_2$ and $\mathcal{L} = \bigl( \mathscr{O}_X,
  \mathscr{O}_{X}(1,0), \mathscr{O}_{X}(0,1) \bigr)$ as in
  Example~\ref{exa:FF1}.  Since $L_1 \otimes L_2 = \mathscr{O}_X(1,1)$
  is very ample and the map
  \[
  H^0(X,L_1) \otimes_\kk H^0(X,L_2) \rTo H^0 \bigl(X,
  \mathscr{O}_X(2,2) \bigr)
  \] 
  is surjective, the complete quiver of sections for $\mathcal{L}$ is
  very ample by Corollary~\ref{cor:veryample}.
\end{example}

\begin{example}
  \label{exa:FF2ample}
  Let $X = \FF_2$ and $\mathcal{L} = \bigl( \mathscr{O}_X,
  \mathscr{O}_{X}(1,0), \mathscr{O}_{X}(0,1), \mathscr{O}_{X}(1,1)
  \bigr)$ as in Example~\ref{exa:FF2}.  Since the line bundle $L_1
  \otimes L_2 \otimes L_3 = \mathscr{O}_X(2,2)$ is very ample on $X$
  and the map $H^0(X,L_1) \otimes_\kk H^0(X,L_2) \otimes_\kk
  H^0(X,L_3) \rTo H^0 \bigl(X, \mathscr{O}_X(2,2) \bigr)$ is
  surjective, Corollary~\ref{cor:veryample} implies that the complete
  quiver of sections for $\mathcal{L}$ is very ample.
\end{example}

To see that every list of basepoint-free line bundles belongs to some
very ample quiver of sections, we prove:

\begin{proposition}
 \label{pro:veryample}
 Let $L_1,\dotsc,L_{r-1}$ be basepoint-free line bundles on $X$.  If
 the subsemigroup of $\Pic(X)$ generate by $L_1, \dotsc, L_{r-1}$
 contains an ample line bundle, then there exists a line bundle $L_r$
 such that the complete quiver of sections for $\mathcal{L} =
 (\mathscr{O}_X, L_1, \dotsc, L_r)$ is very ample.
\end{proposition}
 
\begin{proof}
  By choosing $b_1,\dotsc, b_{r-1} \in \NN$ sufficiently large, we may
  assume that the line bundle $L_r := L_1^{b_1} \otimes \dotsb \otimes
  L_{r-1}^{b_{r-1}}$ is very ample and $\mathscr{O}_X$-regular with
  respect to $L_1, \dotsc, L_{r-1}$; for the multigraded definition of
  regularity see \cites{MaclaganSmith, HSS}.  Let $Q$ be the complete
  quiver of sections for $\mathcal{L} = (\mathscr{O}_X, L_1,\dotsc,
  L_r)$.  Since $L_r$ is very ample and $L_0, \dotsc, L_{r-1}$ are
  basepoint-free, it follows that the line bundle $L :=
  \bigotimes_{i\in Q_0} L_i$ is very ample.  Since Theorem~2.1 in
  \cite{HSS} proves that $H^0(X, L_0) \otimes_\kk \dotsb \otimes_\kk
  H^0(X,L_r) \rTo H^0(X, L)$ is surjective,
  Corollary~\ref{cor:veryample} completes the proof.
\end{proof}

\begin{theorem}
  \label{thm:ample}
  If $Q$ is a very ample quiver of sections, then we can recover the
  line bundles $L_0, \dotsc, L_r$ as the restriction of the
  tautological line bundles on $Y = |\mathcal{L}|$.
\end{theorem} 

\begin{proof}
  If we identify $\Wt(Q)$ with $\ZZ^r$ by choosing the weights
  $(\chi_1 - \chi_0, \dotsc, \chi_r - \chi_0)$ as an ordered basis,
  then the projection map $\ZZ^{Q_0} \rTo \ZZ^r$ induces an isomorphism
  between $G$ and the subgroup $\{ (g_0, \dotsc, g_r) \in
  (\kk^*)^{r+1} : g_0 = 1 \}$ of $(\kk^*)^{Q_0}$.  This isomorphism
  determines a $G$-equivariant vector bundle $\bigoplus_{i \in Q_0}
  \mathscr{O}_{\AA^{Q_1}}$; specifically, the $i$-th component
  corresponds to the $S_Y$-module $S_Y(\chi_i - \chi_0)$, where
  $\bigl( S_Y(\theta') \bigr)_\theta = (S_Y)_{\theta' + \theta}$.  If
  follows that the tautological line bundles on $Y$ are
  $\mathscr{O}_Y, \mathscr{O}_Y(\chi_1 - \chi_0), \dotsc,
  \mathscr{O}_Y(\chi_r - \chi_0)$.  Restricting to $\VV(I_Q)$, we
  obtain a $G$-equivariant vector bundle $\bigoplus_{i \in Q_0}
  \mathscr{O}_{\VV(I_Q)}$ where the $i$-th component corresponds to
  $\bigl( S_Y/ I_Q \bigr)(\chi_i - \chi_0)$.  Since $\varphi_Q \colon
  X \rTo Y$ is a closed embedding, Proposition~\ref{pro:ample} implies
  that $\bigl( (S_Y/I_Q)[y^{- \widehat\sigma}] \bigr)_{[0]} \cong
  \bigl( S_X[x^{- \widehat\sigma}] \bigr)_{[0]}$ for all $\sigma \in
  \Sigma_X(d)$.  Hence, the $S_Y$-module $\bigl( S_Y/ I_Q
  \bigr)(\chi_i - \chi_0)$ corresponds to $\pic(\chi_i - \chi_0) =
  L_i$ on $X \cong \varphi_Q(X)$.
\end{proof}

\begin{remark}
  If we identify $\Wt(Q)$ with $\ZZ^r$ by choosing $(\chi_0 - \chi_1,
  \dotsc, \chi_0 - \chi_r)$ as an ordered basis then the restriction
  of the tautological line bundles would yield the inverse line
  bundles $\mathcal{O}_X, L_1^{-1}, \dotsc, L_r^{-1}$.
\end{remark}

\section{Representations of Bound Quivers} 
\label{sec:boundreps}

This section connects the ideal of relations on a quiver of sections
to the geometry of its multilinear series.  Throughout this section,
let $(Q,R)$ be a complete, bound quiver of sections for $\mathcal{L} =
(\mathscr{O}_X,L_1, \dotsc, L_r)$ of line bundles on $X$ and let $Y =
|\mathcal{L}|$ be the multilinear series of $Q$.

The ideal of relations defines an algebraic subset of $\AA^{Q_1}$.
More precisely, the map sending the path $p = a_1 \dotsb a_\ell$ in
$Q$ to the monomial $y_{a_1} \dotsb y_{a_\ell} \in S_Y$ extends to
give a $\kk$-linear map from $\kk Q$ into $S_Y$.  Let $I_R$ be the
ideal in $S_Y$ generated by the image of $R$ under this map; it is a
binomial ideal because $R$ is spanned by differences $p - p' \in \kk
Q$.  Since these differences satisfy $\tail(p) = \tail(p')$, $\head(p)
= \head(p')$, and $\div(p) = \div(p')$, $I_R$ is homogeneous with
respect to the $\Wt(Q)$-grading on $S_Y$ and is contained in $I_Q$.

The following examples illustrate various possible relations between
$I_R$ and $I_Q$.

\begin{example}
  \label{exa:PP2O2}
  Let $X = \PP^1$ and let $\mathcal{L} = \bigl( \mathscr{O}_X,
  \mathscr{O}_X(2) \bigr)$; Example~\ref{exa:2ptY} describes
  $|\mathcal{L}|$.  There are no paths of length greater than
  $1$ in $Q$, so $R$ and $I_R$ are both the zero ideal.  Since $I_Q$
  is the toric ideal associated to the matrix
  \renewcommand{\arraystretch}{0.8}
  \renewcommand{\arraycolsep}{3pt}
  \[
  \left[ \text{\footnotesize $\begin{array}{rrr}
      -1 & -1 & -1 \\ 
      1 &  1 &  1 \\
      2 &  1 &  0 \\
      0 &  1 &  2
    \end{array}$} \right] \, ,
  \]
  we have $I_Q = ( y_0y_2 - y_1^2)$.  Thus, $\VV(I_Q)$ is a closed
  subvariety of $\VV(I_R) = \AA^{3} = \AA^{Q_1}$.
\end{example}

\begin{example}
  \label{exa:FF1Z}
  Let $X = \FF_1$ and $\mathcal{L} = \bigl( \mathscr{O}_X,
  \mathscr{O}_{X}(1,0), \mathscr{O}_{X}(0,1) \bigr)$.  Since
  Example~\ref{exa:FF1} shows that $R = (0)$ and
  Example~\ref{exa:FF2image} shows that $I_Q = (0)$, it follows that
  $I_R = I_Q = (0)$; in other words, $\VV(I_Q) = \VV(I_R) =
  \AA^4=\AA^{Q_1}$.
\end{example}

\begin{example}
  \label{exa:FF2Z}
  Let $X = \FF_2$ and $\mathcal{L} = \bigl( \mathscr{O}_X,
  \mathscr{O}_{X}(1,0), \mathscr{O}_{X}(0,1), \mathscr{O}_{X}(1,1)
  \bigr)$.  Using the description of $R$ in Example~\ref{exa:FF2}, we
  see that
  \begin{align*}
    I_R &= (y_2y_4-y_1y_5, y_4y_8-y_5y_7, y_2y_6-y_3y_8, y_1y_6 -
    y_3y_7) \\ &= (y_2y_4-y_1y_5, y_1y_6-y_3y_7, y_2y_6-y_3y_8,
    y_2y_7-y_1y_8, y_5y_7-y_4y_8) \cap (y_3, y_4, y_5, y_6) \, .
  \end{align*}
  The description of $I_Q$ given in Example~\ref{exa:FF2image} implies
  that $I_Q$ is a primary component of $I_R$.  Geometrically,
  $\VV(I_Q)$ is the unique component of $\VV(I_R)$ not lying in a
  linear subspace.
\end{example}

Let $W = (W_i, w_a)$ be a representation of $Q$.  For any nontrivial
path $p = a_1 \dotsb a_\ell$, the evaluation of $W$ on $p$ is the
$\kk$-linear map $w_p \colon W_{\tail(p)} \rTo W_{\head(p)}$ defined by
the composition $w_p = w_{a_1} \dotsb w_{a_\ell}$.  This definition
extends to $\kk$-linear combinations of paths with a common head and a
common tail.  A representation of the bound quiver $(Q,R)$ is a
representation $W$ of $Q$ such that $w_p - w_{p'} = 0$ for all $p-p'
\in R$.  Consequently, a point in the representation space $\AA^{Q_1}$
for $Q$ corresponds to a representation for $(Q,R)$ if and only it
lies in the subscheme $\VV(I_R)$.  The category of representations of
the bound quiver $(Q,R)$ of dimension vector $\sum_{i\in Q_0}
\chi_i\in \NN^{Q_0}$ is equivalent to the category of $(\kk
Q/R)$-modules that are isomorphic as $\bigl( \bigoplus_{i\in Q_0} \kk
e_i \bigr)$-modules to $\bigoplus_{i\in Q_0} \kk e_i$.

The ideal $I_R$ is homogeneous with respect to the $\Wt(Q)$-grading on
$S_Y$, so the subscheme $\VV(I_R)$ is $G$-invariant where $G =
\Hom_\ZZ\bigl( \Wt(Q), \kk^* \bigr)$.  The GIT-chamber decomposition
of $\Wt(Q) \otimes_\ZZ \QQ$ arising from the $G$-action on $\VV(I_R)$
coarsens that for the $G$-action on $\mathbb{A}^{Q_1}$; see
Remark~\ref{rem:chamber}.  Let $\Theta$ denote the GIT-chamber arising
from the $G$-action on $\VV(I_R)$ containing $\vartheta = \sum_{i\in
  Q_0} (\chi_i-\chi_0)$.  Proposition~5.3 in \cite{King1} shows that,
for $\theta \in \Theta$, the GIT-quotient
\[
\mathcal{M}_\theta(Q,R) := \VV(I_R) \git G = \Proj \Bigl(
\bigoplus\limits_{k \in \NN} \bigl( \tfrac{S_Y}{I_R} \bigr)_{k
  \theta} \Bigr)
\]
is the fine moduli space for $\theta$-stable representations of
$(Q,R)$. Equivalently, if $\bigoplus_{i \in Q_0} \kk e_i$ denotes the
subalgebra of $\End\bigl( \bigoplus_{i \in Q_0} L_i \bigr)$ generated
by the primitive orthogonal idempotents $e_i$ for $i \in Q_0$, then
Proposition~\ref{pro:algebra} implies that $\mathcal{M}_\theta(Q,R)$
is the fine moduli space of $\theta$-stable $\End\bigl( \bigoplus_{i
  \in Q_0} L_i \bigr)$-modules that are isomorphic as $\bigl(
\bigoplus_{i \in Q_0} \kk e_i \bigr)$-modules to $\bigoplus_{i \in Q_0}
  \kk e_i$.

\begin{theorem}
  \label{thm:finequotient}
  If $Q$ is a very ample quiver of sections, then the following are
  equivalent:
  \begin{enumerate}
  \item the ideal $I_Q$ equals ideal quotient $(I_R : B_Y^\infty)$;
  \item for all $\theta \in \Theta$, the map $\varphi_Q$ induces an
    isomorphism from $X$ to $\mathcal{M}_{\theta}(Q,R)$.
  \end{enumerate}
\end{theorem}

\begin{proof}
  The equivalence of (b) and (c) in Proposition~\ref{pro:smooth}
  implies that the moduli space $\mathcal{M}_\theta(Q,R) = \VV(I_R)
  \git G$ is the geometric quotient of $\VV(I_R) \setminus \VV(B_Y)$
  by the group $G$.  Since $Q$ is a very ample quiver of sections,
  Proposition~\ref{pro:ample} proves that the map $\varphi_Q$ induces
  an isomorphism from $X$ to $\varphi_Q(X)$ and
  Proposition~\ref{pro:ideal} establishes that $\varphi_Q(X)$ is the
  geometric quotient of $\VV(I_Q) \setminus \VV(B_Y)$ by $G$.  Since
  $I_Q$ is prime, the locally closed subscheme $\VV(I_Q) \setminus
  \VV(B_Y)$ equals $\VV(I_R) \setminus \VV(B_Y)$ if and only if we
  have $I_Q = I_R : B_Y^\infty$.
\end{proof}

A list $\mathcal{L} = (\mathscr{O}_X, L_1, \dotsc, L_r)$ of line
bundles on $X$ is \emph{fine} if the complete bound quiver of sections
for $\mathcal{L}$ is very ample and satisfies either of the equivalent
conditions in Theorem~\ref{thm:finequotient}.  The next result
shows that every projective toric variety has many fine lists.

\begin{theorem}
  \label{thm:fine}
  Let $L_1, \dotsc, L_{r-2}$ be basepoint-free line bundles on $X$.
  If the subsemigroup of $\Pic(X)$ generated by $L_1, \dotsc, L_{r-2}$
  contains an ample line bundle, then there exist line bundles
  $L_{r-1}$ and $L_r$ such that the list $\mathcal{L} =
  (\mathscr{O}_X, L_1, \dotsc, L_r)$ is fine.
\end{theorem}

\begin{proof} 
  We divide the proof into three parts.
  
  \begin{proof1}
    This part is similar to Proposition~\ref{pro:ample}.  By choosing
    sufficiently large positive integers $b_1, \dotsc, b_{r-2}$, we
    may assume that the line bundle $L_{r-1} := L_1^{b_1} \otimes
    \dotsb \otimes L_{r-2}^{b_{r-1}}$ is $\mathscr{O}_X$-regular with
    respect to $L_1, \dotsc, L_{r-2}$.  Set $L_r := L_{r-1}^2$.  By
    increasing the $b_i$ if necessary, we may also assume that $L_r$
    is very ample.  Let $Q$ be the complete quiver of sections for
    $\mathcal{L} = (\mathscr{O}_X, L_1, \dotsc, L_r)$.  Since $L_r$ is
    very ample and $L_1, \dotsc, L_{r-1}$ are basepoint-free, it
    follows that $L := \bigotimes_{i \in Q_0} L_i$ is very ample.
    Since Theorem~2.1 in \cite{HSS} implies that $H^0(X,L_1)
    \otimes_\kk \dotsb \otimes_\kk H^0(X,L_r) \rTo H^0(X,L)$ is
    surjective, Corollary~\ref{cor:veryample} implies that $Q$ is very
    ample.
  \end{proof1}

  \begin{proof2}
    By definition, $I_Q$ is the toric ideal associated to the map $\pi
    \colon \NN^{Q_1} \rTo \Wt(Q) \oplus \CDiv(X)$.  It suffices by
    Lemma~12.2 in \cite{Sturmfels} to construct a subset $\mathcal{C}$
    that generates the kernel of $\pi$ as an abelian group and
    satisfies $(y^{v_+} - y^{v_-} : v_+- v_- = v \in \mathcal{C})
    \subseteq I_R$.  Since $L_{r-1}$ is
    $\mathscr{O}_X$\nobreakdash-regular and the $b_i$ are positive,
    Theorem~2.1 in \cite{HSS} shows that $H^0(X,L_{r-1}^{\,} \otimes
    L_i^{-1}) \otimes_\kk H^0(X,L_{r-1}) \rTo H^0(X,L_r^{\,} \otimes
    L_i^{-1})$ is surjective for all $1 \leq i \leq r-1$.  Hence,
    every path in $Q$ from $0$ to $r$ passes through $r-1 \in Q_0$.
    Moreover, because $Q$ is complete, the set $\mathcal{A}$ of the
    arrows from $r-1$ to $r$ in $Q$ corresponds to the set of nonzero
    $T_X$-invariant elements in $H^0(X,L_{r-1})$.  Since the set
    $\mathcal{P}$ of paths from $0$ to $r-1$ in $Q$ are labelled by
    nonzero $T_X$-invariant elements in $H^0(X,L_{r-1})$, there is a
    surjective function $\Psi \colon \mathcal{P} \rTo \mathcal{A}$ such
    that $\div(p) = \div\bigl( \Psi(p) \bigr)$ for all $p \in
    \mathcal{P}$.  For $(a, a', p, p') \in \mathcal{A}^2 \times
    \mathcal{P}^2$, we have $\div(p) + \div(a) = \div(p') + \div(a)$
    if and only if $p a - p' a' \in R$; set 
    \[
    \mathcal{C} := \{ f(p) + f(a) - f(p') - f(a') : pa-p'a' \in R \}
    \subseteq \Cir(Q) \, ,
    \] 
    where $f(\gamma)$ is the element in $\ZZ^{Q_1}$ associated to a
    walk $\gamma$ in $Q$ defined in \S\ref{sub:quivers}.

    To analyze $\mathcal{C}$, we use an ``elongation'' operation on
    circuits in $Q$.  Because $Q$ is acyclic, we may assume every
    circuit $\gamma = \alpha_1^{\,} \alpha_2^{-1} \alpha_3^{\,} \dotsb
    \alpha_{2\ell -1}^{\,} \alpha_{2 \ell}^{-1}$ is an alternating
    sequence of forward paths $\alpha_1, \alpha_3, \dotsc,
    \alpha_{2\ell -1}$ and backward paths $\alpha_2^{-1},
    \alpha_4^{-1}, \dotsc, \alpha_{2 \ell}^{-1}$.  Since $Q$ is
    connected, there exists at least one path from the unique source
    $0$ to each $i \in Q_0$.  Similarly, the choice of $L_r$ implies
    that there is at least one path from each $i \in Q_0$ to $r$.  Let
    $\widehat{\gamma}$ denote a closed walk obtained from the circuit
    $\gamma$ via the following procedure: for $i=1,3, \dotsc,
    2\ell-1$, choose a path $\beta_i$ from $\head(\alpha_i)$ to $r$;
    for $i = 2, 4, \dotsc, 2\ell$, choose a path $\beta_i$ from $0$ to
    $\tail(\alpha_i)$; let $\beta_{0} = \beta_{2 \ell}$ and set
    $\widehat{\gamma} := \beta_{0}^{\,} \alpha_1^{\,}\beta_1^{\,}
    \beta_1^{-1}\alpha_2^{-1} \beta_{2}^{-1}\beta_2^{\,} \alpha_3^{\,}
    \beta_3^{\,} \beta_3^{-1} \dotsb \alpha_{2 \ell}^{-1}
    \beta_{2\ell}^{-1}$.  Observe that $\widehat\gamma =
    \widehat{p}_1^{\,} \widehat{p}_2^{-1} \widehat{p}_3^{\,} \dotsb
    \widehat{p}_{2\ell -1}^{\,} \widehat{p}_{2 \ell}^{-1}$ is an
    alternating sequence of forward and backward paths between $0$ and
    $r$, where $\widehat{p}_i = \beta_{i-1}\alpha_{i}\beta_{i}$ for
    odd $i$, $\widehat{p}_i^{-1} =
    \beta_{i-1}^{-1}\alpha_i^{-1}\beta_{i}^{-1}$ for even $i$.  Thus,
    we have $\div(\gamma) = \div(\widehat\gamma) \in \CDiv(X)$ and
    \[
    f(\gamma) = f(\widehat\gamma) = \sum\limits_{i=0}^{\ell} \bigl(
    f(\widehat{p}_{2i-1}) \bigr) - \sum\limits_{i=0}^{\ell} \bigl(
    f(\widehat{p}_{2i}) \bigr) \in \Cir(Q)
    \] 
    where $f(\widehat{p}_i) \in \NN^{Q_1}$.

    To see that $\mathcal{C}$ spans the lattice $\Ker(\pi)$, fix $u
    \in \Ker(\pi)$.  Since $\pi = (\inc, \div)$, the exact sequence
    \eqref{eq:graphses} implies that $u \in \Cir(Q)$.  Theorem~5.2 in
    \cite{Biggs} shows that the circulation lattice is generated by
    the circuits, so $u = \sum_{i} f(\gamma_i)$ where $\gamma_i$ is a
    circuit in $Q$.  By using the elongation operation and regrouping
    the sum, we have $u = \sum_{i} \bigl( f(\widehat{p}_{2i-1}) -
    f(\widehat{p}_{2i}) \bigr)$ where each $\widehat{p}_i$ is a path
    from $0$ to $r$.  Since each path in $Q$ from $0$ to $r$ passes
    through $r-1$, it follows that $\widehat{p}_i = p_i a_i$ where
    $p_i \in \mathcal{P}$ and $a_i \in \mathcal{A}$.  Hence, we have
    $u = \textstyle\sum\nolimits_{i} \bigl( f(p_{2i-1}) + f(a_{2i-1})
    - f(p_{2i}) - f(a_{2i}) \bigr)$.  We decompose this expression
    into a sum of elements in $\mathcal{C}$ by exploiting properties
    of the sets $\mathcal{A}$ and $\mathcal{P}$ arising from our
    choice of $L_{r-1}$. To be specific, let $I'$ be the toric ideal
    associated to the map $\NN^{\mathcal{A}} \rTo \CDiv(X)$ given by
    $\chi_{a} \mapsto \div(a)$.  The identification of $\mathcal{A}$
    with the nonzero $T_X$-invariant elements in $H^0(X,L_{r-1})$
    implies that $I'$ is the ideal of $\varphi_{|L_{r-1}|}(X)$.  Since
    $L_{r-1}$ is $\mathscr{O}_X$-regular, Theorem~1.1 in \cite{HSS}
    establishes that $I'$ is generated by quadrics.  Hence, if
    $\kk[z_a : a \in \mathcal{A}] = \kk[\NN^{\mathcal{A}}]$ and
    \[
    \mathcal{R} := \{ (a_0,a_1,a_2,a_3) \in \mathcal{A}^4 : \div(a_0)
    + \div(a_1) = \div(a_2) + \div(a_3) \} \, ,
    \] 
    then we have $I' = \, \bigr( z_{a_0}z_{a_1} - z_{a_2}z_{a_3} :
    (a_0,a_1,a_2,a_3) \in \mathcal{R} \bigr)$.  For $1 \leq i \leq
    \ell$, we define $c_i := f(p_{2i-1}) + f(\Psi(p_{2i})) - f(p_{2i})
    - f(\Psi(p_{2i-1}))$.  Since $\div(p_i) = \div \bigl( \Psi(p_i)
    \bigr)$ for $1 \leq i \leq \ell$, each $c_i$ belongs to
    $\mathcal{C}$ and we have
    \begin{align*}
      u &= \sum\limits_i \bigl( f(p_{2i-1}) + f(a_{2i-1}) - f(p_{2i})
      - f(a_{2i}) \bigr) \\ &= \sum\limits_i \bigl( f(\Psi(p_{2i-1}))
      + f(a_{2i-1}) - f(\Psi(p_{2i})) - f(a_{2i}) \bigr) +
      \sum\limits_i c_i \, .
    \end{align*}
    Given $v \in \ZZ^{\mathcal{A}}$ satisfying $\div(v) = 0$,
    Theorem~5.3 in \cite{Sturmfels} applied to the generators of the
    toric ideal $I'$ yields
    \begin{align*}
      v = \!\! \sum\limits_{(a_0,a_1,a_2,a_3) \in \mathcal{R}'} \!\! (
      \chi_{a_0} + \chi_{a_1} - \chi_{a_2} - \chi_{a_3} ) = \!\!
      \sum\limits_{(a_0,a_1,a_2,a_3) \in \mathcal{R}'} \!\! \bigl(
      f(a_0) + f(a_1) - f(a_2) - f(a_3) \bigr) \, ,
    \end{align*}
    where $\mathcal{R}'$ is a multiset of elements from $\mathcal{R}$.
    Applying this to the difference $u - \sum_i c_i$ implies that $u =
    \textstyle\sum\nolimits_{(a_0,a_1,a_2,a_3) \in \mathcal{R}'} \bigl(
    f(a_0) + f(a_1) - f(a_2) - f(a_3) \bigr) + \textstyle\sum\nolimits_i
    c_i$ where $\mathcal{R}'$ is some multiset of elements from
    $\mathcal{R}$.  For $1 \leq i \leq \ell$, set $c_i' :=
    f(\widetilde{a}_0) + f(a_2) - f(\widetilde{a}_2) - f(a_0)$ where
    $\widetilde{a}_j \in \mathcal{P}$ satisfies $\Psi(\widetilde{a}_j)
    = a_j$ for $j = 0, 2$.  Hence, we obtain
    \begin{align*}
      u &= \sum\limits_{(a_0,a_1,a_2,a_3) \in \mathcal{R}'} \bigl(
      f(\widetilde{a}_0) + f(a_1) - f(\widetilde{a}_2) - f(a_3) \bigr)
      - \sum\limits_i c_i' + \sum\limits_i c_i
    \end{align*}
    where $\mathcal{R}'$ is a multiset of elements from $\mathcal{R}$.
    Since each summand belongs to $\mathcal{C}$, we see that $u$ lies
    in the lattice spanned by $\mathcal{C}$ and $I_Q = \bigl( I_R :
    (\prod_{a \in Q_1} y_a)^\infty \bigr)$.
  \end{proof2}
  
  \begin{proof3}
    For $u \in \ZZ^{Q_1}$, set $y^u := \prod_{a \in Q_1} y_a^{u_a}$.
    For any subset $Q_1' \subseteq Q_1$, let $f(Q_1') := \sum_{a \in
      Q_1'} \chi_a$ so that $y^{f(Q_1')} := \prod_{a \in Q_1'} y_a$.
    The definition of $B_Y$ given in Proposition~\ref{pro:smooth}~(b)
    implies that $(I_R : B_Y^\infty) = \bigcap_{Q'} \bigl( I_R :
    (y^{f(Q'_1)} )^\infty \bigr)$ where the intersection is over all
    spanning tree $Q' \subseteq Q$ rooted at $0$.  Because $\bigl( I_R
    : ( y^{f(Q_1')})^\infty \bigr)$ is a subset of $\bigl( I_R : (
    \prod_{a \in Q_1} y_a )^\infty \bigr)$, it is enough to show that
    each spanning tree $Q' \subseteq Q$ rooted at $0$ satisfies
    $\bigl( I_R : ( y^{f(Q_1')} )^\infty \bigr) = \bigl( I_R :
    (\prod_{a \in Q_1} y_a )^\infty \bigr)$.

    The technique of proof follows Example~2.3 in \cite{BLG}.  By
    increasing the $b_i$ from part~1 if necessary, we may assume that
    there exists $s \in H^0(X,L_{r-1})$ such that the corresponding
    lattice point in the polytope $P$ associated $L_{r-1}$ lies in the
    interior.  Fix a spanning tree $Q'$ rooted at $0$.  Let $a_s$ be
    the unique arrow $a_s \in \mathcal{A} \cap Q_1'$ and let $p_s$ be
    a path in $Q$ satisfying $\div(p_s) = \div(s) = \div(a_s)$.
    Observe that $y_{a_s}$ is invertible in $S_Y[y^{-f(Q_1')}]$.  For
    any path $p \in \mathcal{P}$, we have $pa_s - p_s \Psi(p) \in R$,
    so $y^{f(p)} - y^{f(p_s)}y_{\Psi(p)}y_{a_s}^{-1}$ belongs to $I_R
    S_Y[y^{-f(Q_1')}]$.  Hence, for any $p,p' \in \mathcal{P}$ and $a,
    a' \in \mathcal{A}$,
    \begin{align*}
      y^{f(p)}y_a - y^{f(p')}y_{a'} &=
      y^{f(p_s)}y_{\Psi(p)}y_ay_{a_s}^{-1} -
      y^{f(p_s)}y_{\Psi(p')}y_{a'}y_{a_s}^{-1} =
      y^{f(p_s)}y_{a_s}^{-1} ( y_{\Psi(p)} y_a - y_{\Psi(p')} y_{a'} )
    \end{align*}
    in $S_Y[y^{- f(Q_1')}]/ I_R S_Y[y^{-f(Q_1')}]$.  Thus, if 
    \[
    J := I_R + \bigl( y_{a_0}y_{a_1} - y_{a_2}y_{a_3} :
    (a_0,a_1,a_2,a_4) \in \mathcal{R} \bigr) \, ,
    \] 
    then we have $J S_Y[y^{- f(Q_1')}] = I_R S_Y[y^{- f(Q_1')}]$.  

    By assumption, $L_{r-1}$ is ample, so the vertices of the polytope
    $P$ are lattice points.  Since the lattice point corresponding to
    $s$ lies in the interior of $P$, we can express $s$ in the form
    $\tfrac{1}{k}\sum_j c_j s_j$ where $k$, $c_j$ are positive
    integers and the $s_j$ correspond to the vertices of $P$.  It
    follows that $z_{a_s}^k - \prod_j z_{a_{s_j}}^{c_j}$ lies in the
    toric ideal
    \[
    I' = \bigl( z_{a_0}z_{a_1} - z_{a_2}z_{a_3} : (a_0,a_1,a_2,a_3)\in
    \mathcal{R} \big)
    \] 
    introduced in part~2.  Changing variables from $z$'s to $y$'s, we
    see that $y_{a_s}^k - \prod_j y_{a_{s_j}}^{c_j} \in J$.  Since
    $y_{a_s}$ is invertible in $S_Y[y^{-f(Q_1')}]$, it now follows
    that all the $y_{a_{s_j}}$ are invertible in the quotient ring
    $S_Y[y^{- f(Q_1')}] / I_R S_Y[y^{- f(Q_1')}]$.  Moreover, we see
    that $y_a$ is invertible in $S_Y[y^{- f(Q_1')}] / I_R S_Y[y^{-
      f(Q_1')}]$ for all $a \in \mathcal{A}$ because the lattice point
    of $P$ corresponding to $a \in \mathcal{A}$ is a positive rational
    combination of some vertices and each variable $y_{a_{s_j}}$ that
    corresponds to a vertex is invertible in $S_Y[y^{- f(Q_1')}] / I_R
    S_Y[y^{- f(Q_1')}]$.  

    If $p' \in \mathcal{P}$ is the unique path supported on the set
    $Q_1'$, then $y^{f(p')}$ is invertible in $S_Y[y^{-f(Q_1')}] / I_R
    S_Y[y^{- f(Q_1')}]$.  Let $a' \in \mathcal{A}$ is the unique arrow
    satisfying $\div(p') = \div(a')$.  The previous paragraph shows
    that $y_{a'}$ in invertible in $S_Y[y^{-f(Q_1')}] / I_R S_Y[y^{-
      f(Q_1')}]$.  For any path $p \in \mathcal{P}$, we have $pa' - p'
    \Psi(p) \in R$ and the identity $y^{f(p)} =
    y^{f(p')}y_{\Psi(p)}y_{a'}^{-1}$ in $I_R S_Y[y^{- f(Q_1')}]$.
    Each monomial on the right side of this identity is invertible in
    $S_Y[y^{- f(Q_1')}] / I_R S_Y[y^{- f(Q_1')}]$ which implies that
    every variable that divides $y^{f(p)}$ is also invertible.  Since
    the path $p \in \mathcal{P}$ was arbitrary, we conclude that, for
    all $a \in Q_1$, $y_a$ is invertible in $S_Y[y^{- f(Q_1')}] / I_R
    S_Y[y^{- f(Q_1')}]$.  Therefore, Corollary~2.6 in \cite{BLG}
    implies that $\bigl(I_R : ( y^{f(Q_1')} )^\infty \bigr) = \bigl(
    I_R : (\prod_{a \in Q_1} y_a )^\infty \bigr)$.
  \end{proof3}

  \noindent
  Since part~(1) shows that $Q$ is very ample, and combining parts~(2)
  and (3) proves that $Q$ satisfies condition~(a) in
  Theorem~\ref{thm:finequotient}, we conclude that $\mathcal{L}$ is
  fine.
\end{proof}

\begin{proof}[Proof of Theorem~\ref{thm:2}]
  By combining Theorem~\ref{thm:fine} and
  Theorem~\ref{thm:finequotient}, it follows that there are many list
  $\mathcal{L} = (\mathcal{O}_X, L_1, \dotsc, L_r)$ of line bundles on
  $X$ such that the induced morphism $\varphi_{|\mathcal{L}|} \colon X
  \rTo |\mathcal{L}|$ identifies $X$ with the fine moduli space
  $\mathcal{M}_{\vartheta}(Q,R)$.  Theorem~\ref{thm:ample} implies
  that the tautological bundles on $\mathcal{M}_{\vartheta}(Q,R)$
  coincide with the line bundles $\mathcal{O}_X, L_1, \dotsc, L_r$.
\end{proof} 

\begin{remark}
  A priori, $I_Q$ depends on the divisors labelling the arrows in $Q$.
  However, if $\mathcal{L}$ is fine then $I_Q$ depends only on $I_R$
  and hence only the bound quiver $(Q,R)$.
\end{remark}

\begin{example}
  For $X = \PP^1$, Example~\ref{exa:PP2O2} shows that $\mathcal{L} =
  \bigl( \mathscr{O}_X, \mathscr{O}_X(2) \bigr)$ is not fine.  To
  obtain a fine list, we add an appropriate line bundle: $\mathcal{L}'
  := \bigl( \mathscr{O}_X, \mathscr{O}_X(2), \mathscr{O}_X(4) \bigr)$.
  The complete quiver of sections for $\mathcal{L}'$ appears in
  Figure~\ref{fig:simplestfine}.
  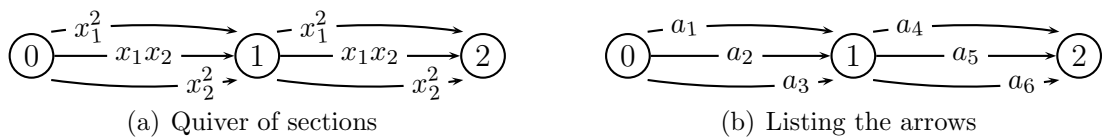
\begin{figure}[!ht]
    \centering
    \mbox{
      \subfigure[Quiver of sections]{
        \psset{unit=1cm}
        \begin{pspicture}(6,1.1) 
          \rput(0,0.5){ 
            \cnodeput(0,0){A}{0}
            \cnodeput(3,0){C}{1} 
            \cnodeput(6,0){E}{2}
            \psset{nodesep=0pt}
            \ncarc[offset=8pt]{->}{A}{C} 
            \ncput*[npos=0.2,nrot=:U]{{\small $x_1^2$}}
            \ncline{->}{A}{C}
            \ncput*[npos=0.5,nrot=:U]{{\small $x_1x_2$}}
            \ncarc[offset=8pt]{<-}{C}{A} 
            \ncput*[npos=0.2,nrot=:180]{{\small $x_2^2$}}
            \ncarc[offset=8pt]{->}{C}{E} 
            \ncput*[npos=0.2,nrot=:U]{{\small $x_1^2$}}
            \ncline{->}{C}{E}
            \ncput*[npos=0.5,nrot=:U]{{\small $x_1x_2$}}
            \ncarc[offset=8pt]{<-}{E}{C} 
            \ncput*[npos=0.2,nrot=:180]{{\small $x_2^2$}}}
        \end{pspicture}}
      \qquad \qquad  
      \subfigure[Listing the arrows]{
        \psset{unit=1cm}
        \begin{pspicture}(6,1.1) 
          \rput(0,0.5){ 
            \cnodeput(0,0){A}{0}
            \cnodeput(3,0){C}{1} 
            \cnodeput(6,0){E}{2}
            \psset{nodesep=0pt}
            \ncarc[offset=8pt]{->}{A}{C} 
            \ncput*[npos=0.2,nrot=:U]{{\small $a_1$}}
            \ncline{->}{A}{C}
            \ncput*[npos=0.5,nrot=:U]{{\small $a_2$}}
            \ncarc[offset=8pt]{<-}{C}{A} 
            \ncput*[npos=0.2,nrot=:180]{{\small $a_3$}}
            \ncarc[offset=8pt]{->}{C}{E} 
            \ncput*[npos=0.2,nrot=:U]{{\small $a_4$}}
            \ncline{->}{C}{E}
            \ncput*[npos=0.5,nrot=:U]{{\small $a_5$}}
            \ncarc[offset=8pt]{<-}{E}{C} 
            \ncput*[npos=0.2,nrot=:180]{{\small $a_6$}}}
        \end{pspicture}}}
    \caption{A fine collection on
      $\mathbb{P}^1_{\kk}$ \label{fig:simplestfine}}
  \end{figure}
  Corollary~\ref{cor:veryample} implies that $|\mathcal{L}'|$ is very
  ample.  Since 
  \begin{align*}
    I_Q &= (y_2^2-y_1y_3,y_5^2-y_4y_6, y_3y_5-y_2y_6, y_2y_5-y_1y_6,
    y_3y_4-y_1y_6,y_2y_4-y_1y_5) \;\; \text{and} \\
    I_R &= (y_3y_5 - y_2y_6, y_3y_4-y_2y_5, y_2y_5-y_1y_6,
    y_2y_4-y_1y_5) 
    = I_Q \cap (y_1,y_2,y_3)\cap (y_4,y_5,y_6) \, ,
  \end{align*}
  it follows that $\mathcal{L}$ is fine.
\end{example}

\begin{example}
  \label{ex:3foldflop}
  Let $X$ be the smooth toric threefold determined by the following
  fan $\Sigma_X$ in $\RR^3$: the rays $\Sigma_X(1)$ are generated by
  the vectors $v_1 := (1,0,0)$, $v_2 := (0,1,0)$, $v_3 := (-1,-1,-1)$,
  $v_4 := (0,1,1)$, $v_5 := (1,0,1)$ and the minimal nonfaces
  correspond to $\{ v_1, v_3, v_4 \}$ and $\{ v_2, v_5 \}$.  The
  induced triangulation of the $2$-sphere is given in
  Figure~\ref{fig:3fold}~(a).  There is a flop
  $\begin{diagram}[inline,size=1.9em,l=>1.9em,tight] X & \rDashto &
    X' \end{diagram}$ where the toric variety $X'$ is the determined
  by the triangulation of $\Sigma_X(1)$ with minimal nonfaces $\{ v_1,
  v_4\}$ and $\{ v_2, v_3, v_5 \}$.  For $(k, \ell) \in \ZZ^2$, write
  $\mathscr{O}_{X}(k,\ell) := \mathscr{O}_{X}(kD_3+\ell D_2) \in
  \Pic(X)$.  The complete quiver of sections for $\mathcal{L} = \bigl(
  \mathscr{O}_{X}, \mathscr{O}_{X}(0,1), \mathscr{O}_{X}(1,0),
  \mathscr{O}_{X}(1,1), \mathscr{O}_{X}(2,0), \mathscr{O}_{X}(2,1)
  \bigr)$ appears in Figure~\ref{fig:3fold}~(b).
  \begin{figure}[!ht]
    \centering
    \mbox{
     \subfigure[Fan]{
        \psset{unit=0.75cm}
        \begin{pspicture}(3,3){
            \cnode*[fillcolor=black](0,0){0pt}{Q1}
            \cnode*[fillcolor=black](3,0){0pt}{Q2}
            \cnode*[fillcolor=black](0,3){0pt}{Q3}
            \cnode*[fillcolor=black](3,3){0pt}{Q4}
            \cnode*[fillcolor=black](0.5,0.5){3pt}{P1}
            \cnode*[fillcolor=black](2.5,0.5){3pt}{P2}
            \cnode*[fillcolor=black](0.5,2.5){3pt}{P5}
            \cnode*[fillcolor=black](2.5,2.5){3pt}{P4}
            \ncline[linewidth=2pt]{-}{P1}{Q1}
            \ncline[linewidth=2pt]{-}{P2}{Q2}
            \ncline[linewidth=2pt]{-}{P5}{Q3}
            \ncline[linewidth=2pt]{-}{P4}{Q4}
            \ncline[linewidth=2pt]{-}{P1}{P2}
            \ncline[linewidth=2pt]{-}{P1}{P5}
            \ncline[linewidth=2pt]{-}{P2}{P4}
            \ncline[linewidth=2pt]{-}{P5}{P4}
            \ncline[linewidth=2pt]{-}{P1}{P4}
            \rput(0,0.7){\psframebox*{$1$}}
            \rput(3,0.7){\psframebox*{$2$}}
            \rput(0,2.3){\psframebox*{$5$}}
            \rput(3,2.3){\psframebox*{$4$}}
            \rput(1.5,3.5){\psframebox*{$3=\infty$}}
          }
      \end{pspicture}}
      \qquad \;\;
      \subfigure[Quiver of sections]{
        \psset{unit=0.75cm}
        \begin{pspicture}(6,3.1)
          \rput(0,0.3){
            \cnodeput(0,0){A}{0}
            \cnodeput(0,2.5){B}{1} 
            \cnodeput(3,0){C}{2}
            \cnodeput(3,2.5){D}{3}
            \cnodeput(6,0){E}{4}
            \cnodeput(6,2.5){F}{5}
            \psset{nodesep=0pt}
            \ncline[offset=5.5pt]{->}{A}{B} \lput*{:180}{$x_2$}
            \ncline[offset=5.5pt]{<-}{B}{A} \lput*{:U}{$x_5$}
            \ncline{->}{A}{C} \lput*{:U}{$x_3$}
            \ncline{->}{B}{D} \lput*{:U}{$x_3$}
            \ncarc[offset=3pt]{->}{B}{C} \lput*{:U}{$x_1$}
            \ncarc[offset=3pt]{<-}{C}{B} \lput*{:180}{$x_4$}
            \ncline[offset=5.5pt]{->}{C}{D} \lput*{:180}{$x_2$}
            \ncline[offset=5.5pt]{<-}{D}{C} \lput*{:U}{$x_5$}
            \ncline{->}{C}{E} \lput*{:U}{$x_3$}
            \ncline{->}{D}{F} \lput*{:U}{$x_3$}
            \ncarc[offset=3pt]{->}{D}{E} \lput*{:U}{$x_1$}
            \ncarc[offset=3pt]{<-}{E}{D} \lput*{:180}{$x_4$}
            \ncline[offset=5.5pt]{->}{E}{F} \lput*{:180}{$x_2$}
            \ncline[offset=5.5pt]{<-}{F}{E} \lput*{:U}{$x_5$}
          }
      \end{pspicture}}
      \qquad \qquad
      \subfigure[Listing the arrows]{
        \psset{unit=0.75cm}
        \begin{pspicture}(6,3.1)
          \rput(0,0.3){
            \cnodeput(0,0){A}{0}
            \cnodeput(0,2.5){B}{1} 
            \cnodeput(3,0){C}{2}
            \cnodeput(3,2.5){D}{3}
            \cnodeput(6,0){E}{4}
            \cnodeput(6,2.5){F}{5}
            \psset{nodesep=0pt}
            \ncline[offset=5.5pt]{->}{A}{B} \lput*{:180}{$a_1$}
            \ncline[offset=5.5pt]{<-}{B}{A} \lput*{:U}{$a_2$}
            \ncline{->}{A}{C} \lput*{:U}{$a_3$}
            \ncline{->}{B}{D} \lput*{:U}{$a_6$}
            \ncarc[offset=3pt]{->}{B}{C} \lput*{:U}{$a_4$}
            \ncarc[offset=3pt]{<-}{C}{B} \lput*{:180}{$a_5$}
            \ncline[offset=5.5pt]{->}{C}{D} \lput*{:180}{$a_8$}
            \ncline[offset=5.5pt]{<-}{D}{C} \lput*{:U}{$a_9$}
            \ncline{->}{C}{E} \lput*{:U}{$a_7$}
            \ncline{->}{D}{F} \lput*{:U}{$a_{12}$}
            \ncarc[offset=3pt]{->}{D}{E} \lput*{:U}{$a_{10}$}
            \ncarc[offset=3pt]{<-}{E}{D} \lput*{:180}{$a_{11}$}
            \ncline[offset=5.5pt]{->}{E}{F} \lput*{:180}{$a_{13}$}
            \ncline[offset=5.5pt]{<-}{F}{E} \lput*{:U}{$a_{14}$}
          }
      \end{pspicture}}
    }
    \caption{Projective threefold admitting a flop \label{fig:3fold}}
  \end{figure}
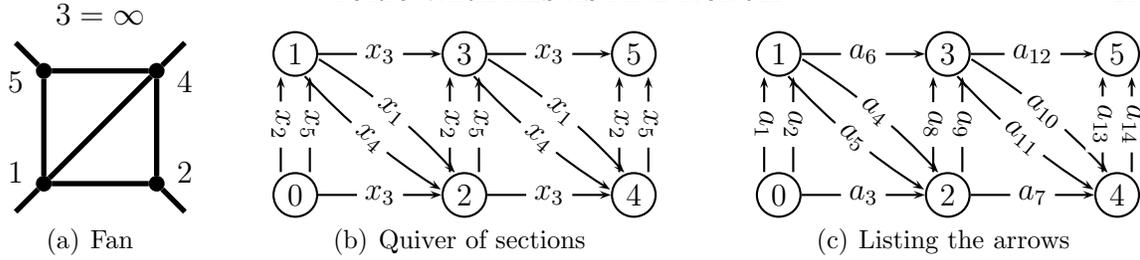

  Since we have
  \begin{align*}
    I_Q &= \left(
      \begin{array}{c}
        y_9y_{13}-y_8y_{14}, y_2y_{13}-y_1y_{14}, y_9y_{12}-y_7y_{14},
        y_8y_{12}-y_7y_{13}, y_5y_{10}-y_4y_{11},\\  y_2y_8-y_1y_9,
        y_5y_7-y_6y_{11}, y_4y_7-y_6y_{10}, y_2y_6-y_3y_9,
        y_1y_6-y_3y_8, \\ y_2y_5y_{12}-y_3y_{11}y_{14}, 
        y_1y_5y_{12}-y_3y_{11}y_{13}, y_2y_4y_{12}-y_3y_{10}y_{14},
        y_1y_4y_{12}-y_3y_{10}y_{13}
      \end{array}
    \right) \\
    I_R &= \left(
      \begin{array}{c} 
        y_9y_{12}-y_7y_{14}, y_8y_{12}-y_7y_{13},
        y_5y_7-y_6y_{11}, y_4y_7-y_6y_{10}, y_2y_6-y_3y_9, \\
        y_1y_6-y_3y_8,   y_9y_{11}y_{13}-y_8y_{11}y_{14},
        y_9y_{10}y_{13}-y_8y_{10}y_{14}, y_5y_9y_{10}-y_4y_9y_{11}, \\
        y_5y_8y_{10}-y_4y_8y_{11}, y_2y_5y_8-y_1y_5y_9,
        y_2y_4y_8-y_1y_4y_9
      \end{array}
    \right) \\
    &= I_Q \cap ( y_{12}, y_{11}, y_{10}, y_7, y_2y_8-y_1y_9,
    y_2y_6-y_3y_9, y_1y_6-y_3y_8 ) \\ &\qquad \; \cap ( y_6, y_5, y_4,
    y_3, y_9 y_{13}-y_8y_{14}, y_9y_{12}-y_7y_{14},y_8y_{12}-y_7y_{13}
    ) \\ &\qquad \; \cap ( y_{12}, y_{11}, y_{10}, y_7, y_6, y_5, y_4,
    y_3 ) \cap ( y_9, y_8, y_7, y_6 )  \, ,
  \end{align*}
  it follows that $\mathcal{L}$ fine.  Observe that $\VV(I_Q)$ has
  four components contained in coordinate hyperplanes of $\AA^{Q_1} =
  \AA^{14}$ and the points in these components correspond to
  representations of disconnected subquivers of $Q$.
\end{example} 

\raggedright

\end{document}